\documentclass[MSNbibl,number,citesort,dvips]{arxbj}
\usepackage{mathbh}

\aid{0}
\volume{20}
\issue{2}
\pubyear{2014}
\firstpage{846}
\lastpage{877}
\doi{10.3150/13-BEJ508} 

\makeatletter
\newcommand{\RMo}{\mathrm{o}}
\newcommand{\RMO}{\mathrm{O}}

\newcommand{\mrmd}{\,\mathrm{d}}
\newcommand{\dd}{\mathrm{d}}

\newremark{example}{Example}[section]
\newproclaim{assumption}[example]{Assumption}
\newtheorem{theorem}[example]{Theorem}
\newtheorem{lemma}[example]{Lemma}
\newremark{remark}[example]{Remark}
\newtheorem{corollary}[example]{Corollary}

\newcommand{\cal}{\mathcal}

\newcommand{\iint}{\int\!\!\int}

\newcommand{\V}{\mathbf{V}}
\newcommand{\N}{\mathbb{N}}
\newcommand{\Z}{\mathbb{Z}}
\newcommand{\R}{\mathbb{R}}
\newcommand{\D}{\mathbb{D}}
\newcommand{\oR}{\overline{\mathbb{R}}}
\newcommand{\pr}{\mathbb{P}}
\newcommand{\ex}{\mathbb{E}}
\newcommand{\vari}{\operatorname{\mathbb{V}ar}}
\newcommand{\covi}{\operatorname{\mathbb{C}ov}}
\newcommand{\eins}{\mathbh{1}}

\newcommand{\rc}{\mathrm{rc}}
\newcommand{\loc}{\mathrm{loc}}
\newcommand{\BV}{\mathbb{BV}}

\makeatother

\begin{document}
\begin{frontmatter}

\title{Continuous mapping approach to the asymptotics of $U$- and
$V$-statistics}
\runtitle{Continuous mapping approach to $U$- and $V$-statistics}

\begin{aug}
\author[1]{\fnms{Eric} \snm{Beutner}\corref{}\thanksref{1}\ead[label=e1]{e.beutner@maastrichtuniversity.nl}} \and
\author[2]{\fnms{Henryk} \snm{Z\"{a}hle}\thanksref{2}\ead[label=e2]{zaehle@math.uni-sb.de}}
\runauthor{E. Beutner and H. Z\"{a}hle} 
\address[1]{Department of Quantitative Economics, Maastricht
University, P.O. Box 616, NL-6200
MD Maastricht, The Netherlands. \printead{e1}}
\address[2]{Department of Mathematics, Saarland University, Postfach
151150, D-66041 Saarbr\"{u}cken, Germany. \printead{e2}}
\end{aug}

\received{\smonth{4} \syear{2012}}
\revised{\smonth{12} \syear{2012}}

%
\begin{abstract}
We derive a new representation for $U$- and $V$-statistics. Using this
representation, the asymptotic distribution of $U$- and $V$-statistics can
be derived by a direct application of the Continuous Mapping theorem.
That novel approach not only encompasses most of the results on the
asymptotic distribution known in literature, but also allows for the
first time a unifying treatment of non-degenerate and degenerate $U$- and
$V$-statistics. Moreover, it yields a new and powerful tool to derive the
asymptotic distribution of very general $U$- and $V$-statistics based on
long-memory sequences. This will be exemplified by several astonishing
examples. In particular, we shall present examples where weak
convergence of $U$- or $V$-statistics occurs at the rate $a_n^3$ and
$a_n^4$, respectively, when $a_n$ is the rate of weak convergence of
the empirical process. We also introduce the notion of asymptotic
(non-) degeneracy which often appears in the presence of long-memory
sequences.
\end{abstract}

%

\end{frontmatter}

\section{Introduction}\label{sectionintroduction}

The study of the asymptotic distribution of $U$- and $V$- (von Mises-)
statistics goes back to Halmos \cite{Halmos1946}, Hoeffding \cite{Hoeffding1948}
and von Mises \cite{vonMises1947}. Different approaches have been
proposed to
obtain the asymptotic distribution of these statistics. The
most-used one is certainly based on the Hoeffding decomposition of a
$U$-statistic; see, for instance,
Dehling \cite{Dehling2006}, Denker \cite{Denker1985},
Koroljuk and Borovskich \cite{KoroljukBorovskich}, Lee \cite{Lee1990},
Serfling \cite{Serfling1980}.
Recently, Beutner and Z{\"a}hle \cite{BeutnerZaehle2011} showed that
the asymptotic
distribution of $U$- and $V$-statistics can be obtained by using the
concept of quasi-Hadamard differentiability introduced in
Beutner and Z{\"a}hle \cite{BeutnerZaehle2010a}. This concept led to
new results for $U$-
and $V$-statistics based on weakly dependent data and was shown in
Beutner, Wu and Z{\"a}hle \cite{BeutnerWuZaehle} to be even suitable
for a certain class of
$U$- and $V$-statistics based on long-memory sequences. However, a
general result that allows to deduce non-central limit theorems for
general $U$- and $V$-statistics based on long-memory sequences is still
missing. This is due to the fact that for long-memory sequences
several parts of the Hoeffding decomposition may contribute to the
limiting distribution; see Dehling and Taqqu \cite{DehlingTaqqu1991}
and our
discussion before Corollary \ref{corollarytothmnclt-1} in
Section \ref{sectionapplication}.

In this article, we derive a new representation of $U$- and
$V$-statistics. Based on this representation, the asymptotic
distribution of $U$- and $V$-statistics, subject to certain regularity
conditions, can be inferred by a direct application of the
Continuous Mapping theorem. It turns out that the continuous mapping
approach does not only cover the majority of the results known in
literature and allows a unifying treatment of non-degenerate and
degenerate $U$- and $V$-statistics (see also Section \ref
{secnotionofdegeneracyandnon-dgeneracy} for the definitions of non-degeneracy
and degeneracy), but also supplements the existing theorems for $U$-
and $V$-statistics based on long-memory sequences. We shall further
see that the continuous mapping approach allows us to establish
strong laws for $U$- and $V$-statistics. Using the continuous mapping
approach it will also be seen that, once the new representation is
established, the asymptotic distributions of several degenerate $U$-
and $V$-statistics that are usually derived on a case by case basis,
are direct consequences of more general results. Finally, we will
demonstrate that, under certain conditions on the kernel, the
continuous mapping approach is also suitable to derive the
asymptotic distribution of two-sample $U$- and $V$-statistics; see Beutner and Z\"{a}hle \cite{Beutnersuppl}.

To explain our approach, we first of all recall that $U$- and
$V$-statistics (of degree
2) are non-parametric estimators for the characteristic
%
\begin{equation}
\label{ucharacteristicdef} V_g(F):=\iint g(x_1,x_2)
\mrmd F(x_1)\mrmd F(x_2)
\end{equation}
of a distribution function (df) $F$ on the real line, where $g\dvtx \R^2\to
\R$ is some measurable function and it is assumed
that the double integral in (\ref{ucharacteristicdef}) exists. Given
a sequence $(X_i)_{i\in\N}$ of random variables on
some probability space $(\Omega,{\cal F},\pr)$ being identically
distributed according to $F$, the $V$-statistic based on $F_n$ is given by
%
\begin{equation}
\label{vstatisticdef} V_g(F_n)=\iint
g(x_1,x_2) \mrmd F_n(x_1)\mrmd F_n(x_2),
\end{equation}
where $F_n$ denotes some estimate of $F$ based on $X_1,\ldots,X_n$,
and it is assumed that the integral in (\ref{vstatisticdef}) exists
for all $n\in\N$. The corresponding $U$-statistic is given by
$U_{g,n}:=\frac{1}{n(n-1)}\sum_{i=1}^n\sum_{j=1,\not
=i}^ng(X_i,X_j)$. 
Assuming\vspace*{1pt} $\iint|g(x_1,x_2)| \mrmd F_n(x_1)\mrmd F(x_2)<\infty$, we obtain from
(\ref{vstatisticdef}) the decomposition
%
\begin{eqnarray}\label{EqHoeffdingdecomposition}
V_g(F_n)-V_g(F) 
& = & \sum_{i=1}^2
\biggl(\int g_{i,F}(x) \mrmd F_n(x)-\int g_{i,F}(x) \mrmd F(x)
\biggr)
\nonumber\\[-8pt]\\[-8pt]
& &{} + \iint g(x_1,x_2) \mrmd (F_n-F)
(x_1)\mrmd (F_n-F) (x_2)\nonumber
\end{eqnarray}
with $g_{1,F}(\cdot):=\int g( \cdot,x_2) \mrmd F(x_2)$ and
$g_{2,F}(\cdot):=\int g(x_1, \cdot) \mrmd F(x_1)$. This decomposition is
sometimes called von Mises decomposition of $V_g(F_n)-V_g(F)$; see
Koroljuk and Borovskich \cite{KoroljukBorovskich}, page 40. If $F_n$ is
the empirical df $\hat
F_n:=\frac{1}{n}\sum_{i=1}^n\eins_{[X_i,\infty)}$ of $X_1,\ldots,X_n$,
then the first line and the second line on the right-hand side of
(\ref{EqHoeffdingdecomposition}) are the linear part and the
degenerate part of the von Mises decomposition, respectively. In this
case, the linear part of the von Mises decomposition coincides with the
linear part of the Hoeffding decomposition of $U_{g,n}-V_g(F)$, and the
degenerate part differs from the degenerate part of the Hoeffding
decomposition of $U_{g,n}-V_g(F)$ only by
%
\begin{equation}
\label{diffdegparthoeffvonmises} \frac{1}{n^2(n-1)}\sum
_{i=1}^n\sum_{j=1,\not=i}^ng(X_i,X_j)-
\frac
{1}{n^2}\sum_{i=1}^ng(X_i,X_i).
\end{equation}
While the linear part can usually be treated by a central limit theorem
(applied to the random variables $Y_i:=g_{1,F}(X_i)+g_{2,F}(X_i)$,
$i\in\N$), it is exactly the degenerate part that causes the main
difficulties in deriving the asymptotic distribution of $U$- and $V$-statistics.

Now let us suppose that we may apply a one-dimensional
integration-by-parts formula to the two summands in the first line on the
right-hand side of (\ref{EqHoeffdingdecomposition}) and a
two-dimensional integration-by-parts formula to the second line on
the right-hand side of (\ref{EqHoeffdingdecomposition}). Notice
that this assumption in particular implies that $g_{1,F}$ and
$g_{2,F}$ generate (possibly signed) measures on $\R$ and that $g$
generates a (possibly signed) measure on $\R^2$. We then have the
following representation (assuming that expressions like
$\lim_{x\to\infty}(F_n-F)(x) g_{i,F}(x)$ are equal to zero
$\pr$-a.s.)
%
\begin{eqnarray}\label{EqBZdecomposition}
V_g(F_n)-V_g(F) & = & -\sum
_{i=1}^2\int\bigl[F_n(x-)-F(x-)\bigr]
\mrmd g_{i,F}(x)
\nonumber\\[-8pt]\\[-8pt]
& &{} + \iint(F_n-F) (x_1-) (F_n-F)
(x_2-) \mrmd g(x_1,x_2),\nonumber
\end{eqnarray}
where we refer to the sum of the first two lines on
the right-hand side as the linear part and to the last line on the
right-hand side as the degenerate part of the representation. Of
course, they coincide with the linear part and the degenerate part
of the von Mises decomposition (\ref{EqHoeffdingdecomposition}).
The representation (\ref{EqBZdecomposition}) is the sum of the
three mappings
%
\begin{eqnarray}\label{eqdefinitionofPhii}
\Phi_{i,g}\dvtx  \V&\longrightarrow&\R,\qquad  \Phi_{i,g}(f):=-\int
f(x-) \mrmd g_{i,F}(x),\qquad i=1,2,
\nonumber\\[-8pt]\\[-8pt]
\Phi_{3,g}\dvtx  \V&\longrightarrow&\R,\qquad \Phi_{3,g}(f):=\iint
f(x_1-)f(x_2-) \mrmd g(x_1,x_2)\nonumber
\end{eqnarray}
applied to $F_n-F$, where $\V$ is some suitable space consisting of
c\`adl\`ag functions on $\oR$. Of course, if $g$ is symmetric then
(\ref{EqBZdecomposition}) can be represented using two mappings
only. Now, on one hand, if the functions $g_{i,F}$, $i=1,2$,
generate \textit{finite} (possibly signed) measures on $\R$, and if
$g$ generates a \textit{finite} (possibly signed) measure on $\R^2$,
then the mappings $\Phi_{i,g}$, $i=1,2,3$, are continuous if we
endow $\V$ with the uniform sup-metric
$d_\infty(f,h):=\|f-h\|_\infty$. On the other hand, if the (possibly
signed) measure generated by $g_{i,F}$ is not finite but only
$\sigma$-finite, then the map $\Phi_{i,g}$ is obviously not
continuous w.r.t. the uniform sup-metric $d_\infty$. However, if we
assume, for example, $\int(1 \slash\phi(x)) |\mrmd g_{i,F}|(x)<\infty$
for $i=1,2$, where $\phi\dvtx \R\to[1,\infty)$ is any continuous
function, and $|\dd g_{i,F}|$ denotes the total variation measure
generated by $g_{i,F}$, then we still have $\Phi_{i,g}(f_n)
\rightarrow\Phi_{i,g}(f)$ for $i=1,2$ when the sequence $(f_n)$
converges to $f$ in the weighted sup-metric
$d_\phi(f,g):=\|(f-h)\phi\|_\infty$. If in addition $\iint
(\phi(x_1)\phi(x_2))^{-1} |\dd g|(x_1,x_2)<\infty$, then we also have
$\Phi_{3,g}(f_n) \rightarrow\Phi_{3,g}(f)$ when the sequence
$(f_n)$ converges to $f$ in the weighted sup-metric $d_\phi$, and
$|\dd g|$ denotes the total variation measure generated by $g$. That
is, under appropriate conditions, we have
%
\begin{equation}
\label{EqBZdecompositionwithPhis} a_n \bigl(V_g(F_n)-V_g(F)
\bigr) = \sum_{i=1}^2 \Phi_{i,g}
\bigl(a_n(F_n-F) \bigr) + \Phi_{3,g} \bigl(
\sqrt{a_n}(F_n-F) \bigr)
\end{equation}
with \textit{continuous} mappings $\Phi_{i,g}$, $i=1,2,3$, and $a_n$ a
strictly positive real number. Therefore, once the representation
(\ref{EqBZdecomposition}) has been established, one only needs
weak convergence of the process $a_n(F_n-F)$ w.r.t. the weighted
sup-metric $d_\phi$ to make use of the Continuous Mapping theorem.
The latter is not problematic. For instance, weak convergence of
empirical processes w.r.t. weighted sup-metrics has been
established under various conditions; see, for instance,
Beutner, Wu and Z{\"a}hle \cite{BeutnerWuZaehle}, Chen and Fan \cite
{ChenFan2006},
Shao and Yu \cite{ShaoYu1996}, Shorack and Wellner \cite
{ShorackWellner1986}, Wu \cite{Wu2003,Wu2008}, Yukich \cite{Yukich1992}.
One of the advantages of this approach lies in the fact that weak
convergence of $a_n(F_n-F)$ implies that $\sqrt{a_n}(F_n-F)$
converges in probability to zero, and hence that
$\Phi_{3,g}(\sqrt{a_n}(F_n-F))$ converges in probability to zero.
Thus, with the continuous mapping approach we can easily deal with
the degenerate part of a non-degenerate $V$-statistic. For a
degenerate $V$-statistic the linear part vanishes, that is,
$\Phi_{1,g}\equiv\Phi_{2,g}\equiv0$, and so in this case (\ref
{EqBZdecompositionwithPhis}) multiplied by $a_n$ reads as
$a_n^2(V_g(F_n)-V_g(F))=\Phi_{3,g}(a_n(F_n-F))$. That is, with the
continuous mapping approach we can also easily deal with the degenerate
part of a degenerate $V$-statistic. Moreover, the continuous mapping
approach can also provide a simple way to derive the asymptotic
distribution of a $V$-statistic when both terms of the von Mises
decomposition contribute to the asymptotic distribution, or when the
scaling sequence has to be chosen as the cube $(a_n^3)$ or the fourth
power $(a_n^4)$ of the scaling sequence $(a_n)$ of $\hat F_n-F$. This
will be illustrated in Section~\ref{sectionapplication} in the context
of long-memory data; see Examples \ref{examplevariance-longmem},
\ref{exampleartificial-longmem} and \ref{exampletfs-longmem}. In this
context the notion of asymptotic degeneracy, to be introduced in
Section \ref{secnotionofdegeneracyandnon-dgeneracy}, plays a crucial
role.

Next,\vspace*{1pt} let us briefly discuss the relation between the asymptotic
distribution of a $U$-statistic $U_{g,n}$ with that of the $V$-statistic
$V_g(\hat F_n)$, where as before $\hat F_n$ denotes the empirical
df. Since the linear parts of the von Mises decomposition and the
Hoeffding decomposition coincide, and the degenerate parts of these
decompositions differ only by the term in (\ref
{diffdegparthoeffvonmises}), it can be shown easily that a non-degenerate
$V$-statistic $V_g(\hat F_n)$ w.r.t. $(g,F)$ and the corresponding
non-degenerate $U$-statistic $U_{g,n}$ have the same asymptotic
distribution if $a_n=\RMo(n)$ and $\ex[|g(X_1,X_1)|]<\infty$; see, for
example, Beutner and Z{\"a}hle \cite{BeutnerZaehle2011}, Remark 2.5.
Hence, in the
non-degenerate case the asymptotic distribution of both $U$-statistics
and $V$-statistics can be derived from (\ref{EqBZdecomposition}). On
the other hand, in the degenerate case they differ by a constant if
a LLN holds for $\sum_{i=1}^ng(X_i,X_i)$. In fact, this follows, because
$n(U_{g,n}-V_g(F))$ equals $n(V_g(F_n)-V_g(F))-\frac{n}{n-1}\frac
{1}{n}\sum_{i=1}^ng(X_i,X_i)
+\frac{n}{n-1}(V_g(F_n)-V_g(F))+\frac{n}{n-1}V_g(F)$.

We stress that for kernels $g$ that are locally of
unbounded variation the asymptotic distribution of the corresponding
$U$-statistic cannot be obtained by the continuous mapping approach;
see Beutner and Z{\"a}hle \cite{Beutnersuppl}, Remark 1.1. However, if
such a kernel is
non-degenerate a non-trivial limiting distribution of the
corresponding $U$-statistic can be obtained either by the Hoeffding
decomposition or by using the approach of Beutner and Z{\"a}hle \cite
{BeutnerZaehle2011}
that is based on the concept of quasi-Hadamard differentiability and
the Modified Functional Delta Method. In case that the kernel is
locally of unbounded variation and degenerate, the approach of
Beutner and Z{\"a}hle \cite{BeutnerZaehle2011} also yields little.
However, then the
traditional approach to degenerate $U$- and $V$-statistics that is
briefly recalled after Example \ref{exampleweaklydependentdata}
may lead to a non-trivial limiting distribution.

The rest of the article is organized as follows. In Section \ref
{secnotionofdegeneracyandnon-dgeneracy}, we will discuss the notion of
(non-) degenerate and asymptotically (non-) degenerate $U$- and
$V$-statistcs. In Section \ref{sectioncontmappingapproach}, we will
first give conditions on the kernel $g$, the df $F$ and the estimator
$F_n$ that ensure that the representation (\ref {EqBZdecomposition})
holds (Section \ref{subsectionmainresult}). Thereafter, we will give
interesting examples for kernels $g$ that satisfy these conditions
(Section \ref {subsectionexamplesforthedecomposition}) and apply the
continuous mapping approach to derive weak and strong limit theorems
for $U$- and $V$-statistics based on weakly dependent data (Sections
\ref{secweaklimittheorems} and~\ref{secstronglimittheorems}).
In Section \ref{sectionapplication}, the whole
strength of our approach will be illustrated by deriving non-central
limit theorems for $U$- and $V$-statistics based on strongly dependent
data. The representation (\ref{EqBZdecomposition}) along with a
new non-central limit theorem for the empirical process of a linear
long-memory process offers a very simple way to deduce such
non-central limit theorems. We will present in particular three
astonishing examples. In Example \ref{examplevariance-longmem}
both terms of the von Mises decomposition of a non-degenerate
$V$-statistic contribute to the asymptotic distribution whatever the
true df $F$ is (in the Gaussian case this example is already known
from Dehling and Taqqu~\cite{DehlingTaqqu1991}, Section 3), and in Examples
\ref{exampleartificial-longmem} and \ref{exampletfs-longmem}
the scaling sequences for degenerate $V$-statistics are given by
$(a_n^3)$ and $(a_n^4)$, respectively (and not as usual by the
square $(a_n^2)$) of the scaling sequence $(a_n)$ of $\hat F_n-F$.
A supplemental article Beutner and Z{\"a}hle \cite{Beutnersuppl}
contains a section that
discusses some extensions and limitations of our approach.


\section{The notions of (non-) degeneracy and asymptotic (non-)
degeneracy}\label{secnotionofdegeneracyandnon-dgeneracy}

In this section, we will recall the notion of (non-) degenerate $U$- and
$V$-statistics, and we shall introduce the notion of asymptotically
(non-) degenerate $U$- and $V$-statistics. We will restrict to the case
where $g$, $F$ and $F_n$ admit the von Mises decomposition (\ref
{EqHoeffdingdecomposition}).

The corresponding $V$-statistic $V_g(F_n)$ will be called \textit{non-degenerate} w.r.t. $(g,F)$, if the linear part of the von Mises
decomposition (i.e., the first line on the
right-hand side in (\ref{EqHoeffdingdecomposition})) does not
vanish. The corresponding $V$-statistic $V_g(F_n)$ will be called \textit{degenerate} w.r.t. $(g,F)$ if the linear part of the von Mises
decomposition vanishes, that is, if $\sum_{i=1}^2\int
g_{i,F} \mrmd (F_n-F)=0$ $\pr$-a.s. for every $n\in\N$. This condition
holds in particular if $g_{1,F}\equiv g_{2,F}\equiv0$, or if $F_n$
is a (random) df and both $g_{1,F}$ and $g_{2,F}$ are constant. If
the linear part of the von Mises decomposition does (not) vanish
when $F_n=\hat F_n$, then we also call the corresponding $U$-statistic
$U_{g,n}$ (non-) degenerate w.r.t. $(g,F)$. Recall that it is very
common, mainly in the i.i.d. set-up, to call a $U$-statistic
degenerate if $\vari[g_{i,F}(X_1)]=0$ for $i=1,2$. Notice that, in
this case, this is in line with the convention used here. Indeed, it
is easily seen that $\vari[g_{i,F}(X_1)]=0$ is equivalent to $\int
g_{i,F} \mrmd (\hat F_n-F)=0$ $\pr$-a.s. if $\hat F_n$ is based on an
i.i.d. sequence. Table \ref{tablefinitedegenerate} displays some
examples for non-degenerate and degenerate $U$- and $V$-statistics.

\begin{table}
\caption{Examples for non-degenerate and degenerate $V$-statistics
w.r.t. $(g,F)$}
\label{tablefinitedegenerate}
\begin{tabular*}{\tablewidth}{@{\extracolsep{\fill}}lll@{}}
\hline
Non-degenerate & Gini's mean difference & Example \ref{examplegini} \\
& Variance & Example \ref{examplevariance} \\
[6pt]
Degenerate & Gini's mean difference & Example \ref
{exampleginidegenrate}\\
& (uniform two-point distribution) & \\
& Variance & Example \ref
{exampledegneratevariaceandsquared}(i)\\
& ($4$th central $=$ squared $2$nd central moment of $F$) & \\
& Cram\'{e}r--von Mises & Example \ref{gof}\\
& test for symmetry & Example \ref{tfs}\\
\hline
\end{tabular*}
\end{table}

To introduce the notion of asymptotically (non-) degenerate $U$- and
$V$-statistics, we let $(a_n)\subset(0,\infty)$ be a scaling sequence
such that $a_n(F_n-F)$ converges in distribution to a non-degenerate
limit. The representation (\ref{EqBZdecompositionwithPhis})
indicates that for every (non-degenerate) $V$-statistic $V_g(F_n)$
(w.r.t. $(g,F)$) only the linear part of the von Mises
decomposition may contribute to the limiting distribution of
$a_n(V_g(F_n)-V_g(F))$. If there is a non-trivial limiting
distribution of the linear part weighted by $a_n$, then we call the
$V$-statistic $V_g(F_n)$ \textit{asymptotically non-degenerate} w.r.t.
$(g,F,(a_n))$, and the analogous terminology is used for
$U$-statistics. Of course, every asymptotically non-degenerate $U$- or
$V$-statistic w.r.t. $(g,F,(a_n))$ must also be non-degenerate
w.r.t. $(g,F)$. However, it might happen that the limiting
distribution of the linear part weighted by $a_n$ vanishes. In this
case, we call the $V$-statistic $V_g(F_n)$ \textit{asymptotically
degenerate} w.r.t. $(g,F,(a_n))$, and again the analogous
terminology is used for $U$-statistics. Of course, every degenerate $U$-
or $V$-statistic w.r.t. $(g,F)$ is also asymptotically degenerate
w.r.t. $(g,F,(a_n))$.

For an asymptotically degenerate $U$- or $V$-statistic w.r.t.
$(g,F,(a_n))$ a non-trivial asymptotic distribution can typically be
obtained by weighting the empirical difference by $a_n^2$ instead
of~$a_n$, that is, by considering the limiting distribution of
$a_n^2(V_g(F_n)-V_g(F))$. In this context, two different things may
occur:
\begin{itemize}[(2)]
\item[(1)] The asymptotic distribution of $a_n^2(V_g(F_n)-V_g(F))$ is
non-trivial. In this case, we say that the asymptotically degenerate $U$-
or $V$-statistic w.r.t. $(g,F,(a_n))$ is of \textit{type 1}.
\item[(2)] The asymptotic distribution of $a_n^2(V_g(F_n)-V_g(F))$
is \textit{still} degenerate. In this case, we say that the asymptotically
degenerate $U$- or $V$-statistic w.r.t. $(g,F,(a_n))$ is of \textit{type 2}.
\end{itemize}
It seems that behavior (2) only appears in the presence of
long-memory sequences; for examples, see Examples \ref
{exampleartificial-longmem} and \ref{exampletfs-longmem}. It is worth
pointing out that for an (asymptotically) degenerate $U$- and
$V$-statistics of type 2 a non-trivial limiting distribution can
sometimes be obtained by considering the limiting distribution of
$a_n^p(V_g(F_n)-V_g(F))$ for some $p>2$; see again Examples
\ref{exampleartificial-longmem} and~\ref{exampletfs-longmem}. In case
(1) we can distinguish between the following three
cases:
\begin{itemize}[(1.c)]
\item[(1.a)] Only the degenerate part of the von Mises decomposition
contributes to the limiting distribution of
$a_n^2(V_g(F_n)-V_g(F))$. This is in particular the case, if the $U$- or
$V$-statistic is even degenerate w.r.t. $(g,F)$.
\item[(1.b)] Only the linear part contributes to the limiting
distribution of
$a_n^2(V_g(F_n)-V_g(F))$. This can happen only if the $U$- or $V$-statistic
is asymptotically degenerate w.r.t.
$(g,F,(a_n))$, but non-degenerate w.r.t. $(g,F)$.
\item[(1.c)] Both the linear and the degenerate part contribute to the
limiting distribution of
$a_n^2(V_g(F_n)-V_g(F))$. Again, this can occur only if the $U$- or
$V$-statistic is asymptotically degenerate w.r.t. $(g,F,(a_n))$, but
non-degenerate w.r.t. $(g,F)$.
\end{itemize}
%

In the original version of the manuscript, we guessed that the cases
(1.b) and (1.c) only appear for $U$- and $V$-statistics based on
long-memory sequences. However, a referee provided us with the
following example that shows that behavior (1.c) also occurs for
$m$-dependent sequences.

\begin{example}\label{examplereferee}
Let $(\xi_n)$ and $(\delta_n)$ be independent i.i.d. sequences with
$\pr[\xi_i=0]=\pr[\xi_i=1]=1 \slash2$ and
$\pr[\delta_i=0]=\pr[\delta_i=1]=1 \slash2$. Define the
$1$-dependent sequence $(Z_i)$ by $Z_i:=\xi_i-\xi_{i-1}$ and the
$1$-dependent sequence $(X_i)$ by
\[
X_i:= \cases{ \sqrt{2}, &\quad $Z_i=1$ and $
\delta_i=1$,
\cr
1, &\quad $Z_i=0$ and $\delta_i=1$,
\cr
0, &\quad $Z_i=-1$,
\cr
-1, &\quad $Z_i=0$ and $
\delta_i=0$,
\cr
-\sqrt{2}, &\quad $Z_i=1$ and $
\delta_i=0$.}
\]
We then have $\mu:=\ex[X_i]=0$, $\sigma^2:=\vari[X_i]=1$, and
$(X_i-\mu)^2-\sigma^2=Z_i$. Denote by $\hat{\sigma}_n^2$ the sample
variance, which is the $V$-statistic $V_g(\hat F_n)$ with kernel
$g(x_1,x_2)=\frac{1}{2}(x_1-x_2)^2$. The Hoeffding decomposition of
$\hat{\sigma}_n^2-\sigma^2$ equals
\[
\frac{1}{n}\sum_{i=1}^n
\bigl((X_i-\mu)^2-\sigma^2\bigr)+\Biggl(
\frac{1}{n}\sum_{i=1}^n(X_i-
\mu)\Biggr)^2 = \frac{1}{n} (\xi_n-
\xi_0) + \Biggl(\frac{1}{n} \sum_{i=1}^n
X_i\Biggr)^2,
\]
where the latter identity follows from $(X_i-\mu)^2-\sigma^2=Z_i$,
$\sum_{i=1}^n Z_i=\xi_n-\xi_0$, and $\mu=0$. Thus, we have
$\sqrt{n}(\hat{\sigma}_n^2-\sigma^2)=\frac{1}{\sqrt{n}}(\xi
_n-\xi_0)+(\frac{1}{n^{1/4}} \frac{1}{\sqrt{n}}\sum_{i=1}^n X_i)^2$,
and so we obtain from the central limit theorem along with Slutzky's
lemma, and the fact that $\xi_n-\xi_0$ has the same distribution for
every $n\in\N$, that $\hat\sigma_n^2$ is non-degenerate w.r.t.
$(g,F)$, but asymptotically degenerate w.r.t. $(g,F,(\sqrt{n}))$,
where $F$ refers to the df of the $X_i$. On the other hand, we have $
n(\hat{\sigma}_n^2-\sigma)=\xi_n-\xi_0+(\frac{1}{n}\sum_{i=1}^n X_i)^2$.
As already mentioned $\xi_n-\xi_0$ has the same distribution for
every $n\in\N$, and
$(\frac{1}{\sqrt{n}}\sum_{i=1}^n X_i)^2$ converges in distribution
to a $\chi^2$ distribution with one degree of freedom. Since
$\xi_n-\xi_0$ is the linear part and
$(\frac{1}{\sqrt{n}}\sum_{i=1}^n X_i)^2$ the degenerate part of the
Hoeffding decomposition, the example shows that even for
$m$-dependent sequences both terms of the Hoeffding decomposition
may contribute to the limiting distribution.
\end{example}

Table \ref{tableasymptoticallydegenerate} displays some examples for
asymptotically non-degenerate and degenerate $U$- and $V$-statistics for a
linear process with long-memory. It is worth mentioning that there is a
%
\begin{table}
\caption{Examples for asymptotically non-degenerate and asymptotically
degenerate $U$- and $V$-statistics w.r.t. $(g,F,(a_n))$ for
$a_n:=n^{p(\beta-1/2)}\ell(n)^{-p}$ with $p(2\beta-1)<1$, where the
observations are drawn from a linear process $X_t:=\sum_{s=0}^\infty
a_s\varepsilon_{t-s}$ with $a_s=s^{-\beta}\ell(s)$ for some $\beta
\in(\frac{1}{2},1)$ and some slowly varying $\ell$ (long-memory)}
\label{tableasymptoticallydegenerate}
\begin{tabular*}{\tablewidth}{@{\extracolsep{\fill}}llll@{}}
\hline
Asymptotically & Gini's mean difference & ($p=1$) & Disc. before
Corollary \ref{corollarytothmnclt-1} \\
non-degenerate & & & \\
[6pt]
Asymptotically & (1.a) Cram\'{e}r--von Mises & ($p=2$) & Disc. before
Corollary \ref{corollarytothmnclt-1}\\
degenerate -- type 1 & (1.b) Squared absolute mean & ($p=2$) & Example
\ref{examplesquaredabsoluetvalue-longmem}\\
& (1.c) Variance & ($p=2$) & Example \ref{examplevariance-longmem}\\
[6pt]
Asymptotically & Some artificial kernel & ($p = 3$) & Example \ref
{exampleartificial-longmem}\\
degenerate -- type 2 & test for symmetry & ($p = 4$) & Example \ref
{exampletfs-longmem}\\
\hline
\end{tabular*}
\end{table}
difference between the asymptotic degeneracy of the variance in the
case of a linear long-memory sequence and the $m$-dependent sequence of
Example \ref{examplereferee}. For a long-memory sequence, the
variance is asymptotically degenerate whatever the underlying df is,
whereas for an $m$-dependent sequence it does depend on the underlying
df whether the variance is asymptotically degenerate or not.



\section{\texorpdfstring{Representation (\protect\ref{EqBZdecomposition}): Conditions, examples, and applications}
{Representation (5): Conditions, examples, and applications}}\label{sectioncontmappingapproach}

Let $\D$ be the space of all bounded c\`adl\`ag functions on $\oR$.
Any metric subspace $(\V,d)$ of $\D$ will be equipped with the
$\sigma$-algebra ${\cal V}:={\cal D}\cap\V$ to make it a measurable
space, where ${\cal D}$ is the $\sigma$-algebra generated by the
usual coordinate projections $\pi_x\dvtx \D\to\R$, $x\in\oR$. The roles
of $\V$ and $d$ will often be played by the space $\D_\phi$ of all
$f\in\D$ with $\|f\phi\|_\infty<\infty$ and the weighted sup-metric
$d_\phi(f,h):=\|(f-h)\phi\|_\infty$, respectively, where
$\phi\dvtx \oR\to[1,\infty]$ is any continuous function being real-valued
on $\R$ (henceforth called weight function) and where we use the
convention $0 \cdot\infty:=0$. We will frequently work with the
particular weight function $\phi_\lambda(x):=(1+|x|)^\lambda$ for
fixed $\lambda$.

Further, let $\BV_{\loc,\rc}$ be the space of all functions on $\R$
that are right-continuous and locally of bounded variation, and notice
that every function in $\BV_{\loc,\rc}$ has also left-hand limits. For
$\psi\in\BV_{\loc,\rc}$, we denote by $d\psi^+$ and $d\psi^-$ the
unique positive Radon measures induced by the Jordan decomposition of
$\psi$ and we set $|d\psi|:=d\psi^++d\psi^-$. Analogously, let
$\BV_{\loc,\rc}^2$ be the space of all functions on $\R^2$ that are
upper right-continuous and locally of bounded variation, and for
$\tau\in\BV_{\loc,\rc}^2$, $d\tau^+$, $d\tau^-$ and $|d\tau|$ are
defined analogously to $d\psi^+$, $d\psi^-$ and $|d\psi|$; for details
see the discussion subsequent to Remark \ref
{remarkconditionsintheonedimensioanlintbyparts} below. We shall
interpret integrals as
being over the open intervals $(-\infty,\infty)$ and
$(-\infty,\infty)^2$, that is, $\int=\int_{(-\infty,\infty)}$\vspace*{1pt} and
$\iint=\iint_{(-\infty,\infty)^2}$. Moreover, for a measurable function
$f$ we shall say that the integral of $f$ w.r.t. a signed measure
$\mu$ exists if the four integrals $\int f^+ \mrmd \mu^+$, $\int
f^- \mrmd \mu^+$, $\int f^+ \mrmd \mu^-$ and $\int f^- \mrmd \mu^-$ are all finite,
where $f^+$ and $f^-$ denote the positive and the negative part of $f$,
and $\mu^+$ and $\mu^-$ denote the positive and the negative part of
$\mu$. We denote by $\stackrel{\mathsf d}{\longrightarrow}$ convergence in
distribution in the sense of Pollard~\cite{Pollard1984}, and the Borel
$\sigma$-algebra in $\R$ is denoted by ${\cal B}(\R)$.


\subsection{\texorpdfstring{Conditions for the representation (\protect\ref{EqBZdecomposition})}
{Conditions for the representation (5)}}\label{subsectionmainresult} 

In this section, we provide conditions on $g$, $F$ and the estimate
$F_n$ of $F$ under which the representation (\ref{EqBZdecomposition})
holds true. First of all, we impose assumptions on $g$,
$F$ and $F_n$ that ensure that $V_g(F)$ and $V_g(F_n)$ are well defined.

\begin{assumption}\label{assumptionbasiconUFandFn}
The integral in (\ref{ucharacteristicdef}) exists, the estimate
$F_n$ of $F$ is a non-decreasing c\`adl\`ag process with variation
bounded by 1, and for all $n \in\N$ we have that $\pr$-a.s. $\iint
|g(x_1,\allowbreak x_2)| \mrmd F_n(x_1)\mrmd F_n(x_2)<\infty$.
\end{assumption}

A further minimum requirement for the representation (\ref
{EqBZdecomposition}) is the following.
%
\begin{assumption}\label{assumptionbasiconUFandFn-2}
For all $n \in\N$ we have that $\pr$-a.s. $\iint|g(x_1,x_2)|
\mrmd F_n(x_1)\mrmd F(x_2)<\infty$ and $\iint|g(x_1,x_2)|
\mrmd F(x_1)\mrmd F_n(x_2)<\infty$.
\end{assumption}

Notice that the conditions on $F_n$ imposed by Assumptions \ref
{assumptionbasiconUFandFn}--\ref{assumptionbasiconUFandFn-2} are always
fulfilled if the integral in (\ref{ucharacteristicdef}) exists and
$F_n$ is the empirical df $\hat{F}_n$.
From (\ref{EqHoeffdingdecomposition}) and Lemmas \ref
{Lemmaintegrationbypartsonedimension} and \ref
{Lemmaintegrationbypartstwodimensions} below, we immediately obtain the
following theorem.

\begin{theorem}\label{EqBZdecomposition-valid}
If the assumptions of Lemmas \ref{Lemmaintegrationbypartsonedimension}
and \ref{Lemmaintegrationbypartstwodimensions} (below)
are fulfilled, then the representation (\ref{EqBZdecomposition}) of
$V_g(F_n)-V_g(F)$ holds true $\pr$-a.s. for every $n\in\N$.
\end{theorem}

The following lemma gives conditions that allow to apply almost surely
an integration-by-parts formula (see Beutner and Z{\"a}hle \cite
{BeutnerZaehle2011}, Lemma
B.1) to the two summands in the first line on the right-hand side in
(\ref{EqHoeffdingdecomposition}).

\begin{lemma}\label{Lemmaintegrationbypartsonedimension}
Suppose that:
\begin{itemize}[(a)]
\item[(a)] Assumptions \ref{assumptionbasiconUFandFn}--\ref
{assumptionbasiconUFandFn-2} hold,
\item[(b)] $g_{i,F}\in\BV_{\loc,\rc}$,
\item[(c)] $\int|F_n(x-)-F(x-)| |\dd g_{i,F}|(x) < \infty$ $\pr
$-a.s., for all $n\in\N$ and $i=1,2$,
\item[(d)] $\lim_{|x|\to\infty}(F_n-F)(x) g_{i,F}(x)=0$ $\pr
$-a.s., for all $n\in\N$ and $i=1,2$.
\end{itemize}
Then $\pr$-a.s., for every $n\in\N$,
\begin{eqnarray*}
& & \iint g(x_1,x_2)\mrmd F_n(x_1)\mrmd F(x_2)-
\iint g(x_1,x_2) \mrmd F(x_1)\mrmd F(x_2)
\\
&&\quad = - \int\bigl[F_n(x_1-)-F(x_1-)\bigr]
\mrmd g_{1,F}(x_1),
\\
&& \iint g(x_1,x_2)\mrmd F(x_1)\mrmd F_n(x_2)-
\iint g(x_1,x_2) \mrmd F(x_1)\mrmd F(x_2)
\\
&&\quad = - \int\bigl[F_n(x_2-)-F(x_2-)\bigr]
\mrmd g_{2,F}(x_2).
\end{eqnarray*}
\end{lemma}

\begin{pf}
We only prove the first equation. From Assumptions \ref
{assumptionbasiconUFandFn} and \ref{assumptionbasiconUFandFn-2},
and using Fubini's theorem, we
have that the integrals $\int|g_{1,F}(x_1)| \mrmd F_n(x_1)$ and $\int
|g_{1,F}(x_1)| \mrmd F(x_1)$ are finite, and so $\int|g_{1,F}(x_1)|
|\dd(F_n-F)|(x_1)$ exists. Moreover, for every $n \in\N$, we obtain by
using Fubini's theorem
$\iint g(x_1,x_2)\mrmd F_n(x_1)\mrmd F(x_2)-\iint g(x_1,x_2) \mrmd F(x_1)\mrmd F(x_2)=\int
g_{1,F}(x_1) \mrmd (F_n-F)(x_1)$. Since by assumption (c), we further have that
$\int|F_n(x_1-)-F(x_1-)| |\mrmd g_{1,F}|(x_1)$ exists, the conditions of
Lemma B.1 in Beutner and Z{\"a}hle \cite{BeutnerZaehle2011} are
fulfilled and the result
follows.
\end{pf}

\begin{remark}\label{remarkconditionsintheonedimensioanlintbyparts}
(i) For $F_n=\hat{F}_n$ (recall that $\hat{F}_n$ denotes the
empirical df) conditions (c) and (d) of Lemma \ref
{Lemmaintegrationbypartsonedimension} boil down to conditions on the
tails of $F$ and
$g_{i,F}$, $i=1,2$.

\mbox{}\hphantom{i}(ii) More generally, if for $\pr$-almost every $\omega$ there exist
real numbers $x_{\ell}(\omega)$ and $x_u(\omega)$ such that
$F_n(\omega,x)-F(x)=-F(x)$ for all $x
\leq x_{\ell}(\omega)$, and $F_n(\omega,x)-F(x)=1-F(x)$ for all $x
\geq x_u(\omega)$, then again conditions (c) and (d) of Lemma \ref
{Lemmaintegrationbypartsonedimension} boil
down to conditions on the tails of $F$ and $g_{i,F}$, $i=1,2$.

(iii) Condition (c) holds if $\int \dd g_{i,F}$ exists
for $i=1,2$, and under the conditions of part (ii) of this remark we
have that condition (d) holds if $\|g_{i,F}\|_{\infty} < \infty$ for
$i=1,2$.

(iv) If for some weight function $\phi$ we have that
$d_{\phi}(F_n,F)$ is $\pr$-a.s. finite for every $n\in\N$, and that
$\int1/\phi|\dd g_{i,F}|<\infty$ and
$\lim_{|x|\to\infty}g_{i,F}(x)/\phi(x)=0$ for $i=1,2$, then again
conditions (c) and (d) of Lemma \ref
{Lemmaintegrationbypartsonedimension} hold. Notice also that under the
conditions of part (ii)
of this remark, the condition that $d_{\phi}(F_n,F)$ is $\pr$-a.s.
finite for all $n\in\N$ is a condition on the tails of $F$.

\mbox{}\hphantom{i}(v) If there are $x_{\ell,i}<x_{u,i}$ such that $|g_{i,F}|$ is
non-increasing on $(-\infty,x_{\ell,i}]$ and non-decreasing on
$[x_{u,i},\infty)$ for $i=1,2$, then part (d) of Lemma \ref
{Lemmaintegrationbypartsonedimension} is already implied by
Assumptions \ref{assumptionbasiconUFandFn} and \ref
{assumptionbasiconUFandFn-2}. Indeed, under these
assumptions the integrals $\int|g_{1,F}(x_1)| \mrmd F_n(x_1)$ and $\int
|g_{1,F}(x_1)| \mrmd F(x_1)$
exist, and we have for $x \ge x_{u,i}$ that $
|g_{i,F}(x)(F_n(x)-F(x))|=|{\int_x^{\infty}g_{i,F}(x) \mrmd F(t)}-\int
_x^{\infty}g_{i,F}(x) \mrmd F_n(t)|\leq\int_x^{\infty}|g_{i,F}(t)|
\mrmd F(t)+\int_x^{\infty}|g_{i,F}(t)| \mrmd F_n(t)$,
$i=1,2$. Thus, $\lim_{x\to\infty}(F_n-F)(x)g_{i,F}(x)=0$ for
$i=1,2$. Analogously, we obtain
$\lim_{x\to-\infty}(F_n-F)(x)g_{i,F}(x)=0$. 
\end{remark}

For the functional $\Phi_{3,g}$ to be well defined in the
Lebesgue--Stieltjes sense the kernel $g$ must be upper
right-continuous and locally of bounded variation. For later use and
the reader's convenience we recall the definition of locally bounded
variation. For any function $\tau\dvtx \R^2\to\R$, set $
\mu_\tau({\cal R}_{(x_1,x_2),(y_1,y_2)}):= \tau(x_2,y_2)-\tau
(x_1,y_2)-\tau(x_2,y_1)+\tau(x_{1},y_{1})
$
for every half-open rectangle $\mathcal{R}_{(x_1,x_2),(y_1,y_2)}=
(x_1,x_2] \times(y_1,y_2]$ with $(x_1,x_2) \in\R^2$,
$ x_1 < x_2$, and $(y_1,y_2) \in\R^2$, $
y_1 < y_2$. For a fixed half-open rectangle
$\mathcal{R}=\mathcal{R}_{(a_1,a_2),(b_1,b_2)}= (a_1,a_2] \times
(b_1,b_2]$ in $\R^2$, a pair $P$ of finite sequences
$(x_k)_{k=0,\ldots,n}$ and $(y_{\ell})_{\ell=0,\ldots,m}$ is called
a grid for $\mathcal{R}$ if $a_1=x_0 \leq x_1 \leq\cdots\leq
x_n=a_2$ and $b_1=y_0 \leq y_1 \leq\cdots\leq y_m=b_2$. For any
grid $P$, let
\begin{eqnarray*}
{\cal V}(P,\tau) :\!&= & \sum_{i=1}^n\sum
_{j=1}^m \bigl|\mu_\tau({\cal
R}_{(x_{i-1},x_i),(y_{i-1},y_i)}) \bigr|
\\
& = & \sum_{i=1}^n\sum
_{j=1}^m \bigl|\tau(x_i,y_j)-
\tau(x_{i-1},y_j)-\tau(x_i,y_{j-1})+
\tau(x_{i-1},y_{j-1}) \bigr|.
\end{eqnarray*}

Moreover, let ${\cal V}_\tau(\mathcal{R}):=\sup_{P \in\mathcal{P}}
{\cal V}(P,\tau)$, where $\mathcal{P}$ is the set of all grids for
$\mathcal{R}$. The function $\tau$ is said to be locally of bounded
variation if for every bounded half-open rectangle $\mathcal{R}\subset
\R^2$ we have ${\cal V}_\tau(\mathcal{R})<\infty$, and $\tau$ is said
to be of bounded total variation if there is a constant \mbox{$C>0$} such that
${\cal V}_\tau(\mathcal{R})\le C$ for all bounded half-open rectangles
$\mathcal{R} \subset\R^2$. As mentioned earlier, $\BV_{\loc,\rc}^2$
denotes the space of all upper right-continuous functions
$\tau\dvtx \R^2\to\R$ that are locally of bounded variation, and we use the
two-dimensional Jordan decomposition (see, e.g.,
Ghorpade and Limaye \cite{Inder}, Proposition 1.17) to define $d\tau
^+$, $d\tau^-$ and
$|d\tau|$ similar as $d\psi^+$, $d\psi^-$ and $|d\psi|$. We can now
state the two-dimensional integration-by-parts lemma, which can almost
surely be applied to the second line on the right-hand side in (\ref
{EqHoeffdingdecomposition}).

\begin{lemma}\label{Lemmaintegrationbypartstwodimensions}
Suppose that:
\begin{itemize}[(a)]
\item[(a)] Assumption \ref{assumptionbasiconUFandFn} holds,
\item[(b)] $g \in\BV_{\loc,\rc}^2$, and the functions
$g_{x_1}(\cdot):=g(x_1,\cdot)$ and $g_{x_2}(\cdot):=g(\cdot,x_2)$
are locally of bounded variation for every fixed $x_1$ and $x_2$, respectively,
\item[(c)] $\iint|F_n(x_1-)-F(x_1-)| |F_n(x_2-)-F(x_2-)|
|\dd g|(x_1,x_2) < \infty$ $\pr$-a.s., for all $n\in\N$,
\item[(d)] the following limits exist and equal zero $\pr$-a.s. for
all $n\in\N$:
\begin{eqnarray*}
&&
\lim_{a_1,b_1\to-\infty, a_2,b_2 \to\infty} \biggl[(F_n-F) (b_2)
\int_{a_1}^{a_2} (F_n-F)
(x_1-) \mrmd g_{b_2}(x_1)
\\
&&\hspace*{79.1pt}{} - (F_n-F) (b_1) \int_{a_1}^{a_2}(F_n-F)
(x_1-) \mrmd g_{b_1}(x_1) \biggr],
\\
&&
\lim_{a_1,b_1\to-\infty, a_2,b_2 \to\infty} \biggl[(F_n-F) (a_2)
\int_{b_1}^{b_2} (F_n-F)
(x_2-) \mrmd g_{a_2}(x_2)
\\
&&\hspace*{78.7pt}{} -(F_n-F) (a_1) \int_{b_1}^{b_2}
(F_n-F) (x_2-) \mrmd g_{a_1}(x_2)
\biggr],
\\
&&
\lim_{a_1,b_1\to-\infty, a_2,b_2 \to\infty}  \bigl[(F_n-F) (a_2)
(F_n-F) (b_2) g(a_2,b_2)
\\
&&\hspace*{78.2pt}{} -(F_n-F) (a_1) (F_n-F) (b_2)
g(a_1,b_2)
\\
&&\hspace*{78.2pt}{} -(F_n-F) (a_2) (F_n-F) (b_1)
g(a_2,b_1)
\\
&&\hspace*{78.2pt}{} +(F_n-F) (a_1) (F_n-F) (b_1)
g(a_1,b_1) \bigr].
\end{eqnarray*}
\end{itemize}
Then $\pr$-a.s., for every $n\in\N$,
\[
\iint g(x_1,x_2) \mrmd (F_n-F)
(x_1)\mrmd (F_n-F) (x_2) = \iint(F_n-F)
(x_1-) (F_n-F) (x_2-) \mrmd g(x_1,x_2).
\]
%
\end{lemma}

In part (d) of the lemma the expression ``$ \lim_{a_1,b_1\to-\infty,
a_2,b_2 \to\infty}(\ldots)$'' is understood as convergence of a net
$(\ldots)_{(n_1,n_2,n_3,n_4)\in\N^4}$, with $(-a_1,a_2,-b_1,b_2)$
playing the role of $(n_1,n_2,n_3,n_4)$, where as usual $\N^4$ is
regarded as a directed set w.r.t. the relation $\triangleleft$, and
$(m_1,m_2,m_3,m_4)\triangleleft(n_1,n_2,n_3,n_4)$ means that $m_i\leq
n_i$ for $i=1,\ldots,4$. The analogous interpretations are used for the
other limits.

\begin{pf*}{Proof of Lemma \ref{Lemmaintegrationbypartstwodimensions}}
For two functions $f, h \in\BV_{\loc,\rc}^2$ and every fixed
rectangle $(a_1,a_2] \times(b_1,b_2]$, we have
\begin{eqnarray*}
\int_{a_1}^{a_2} \int_{b_1}^{b_2}
f(x_1,x_2) \mrmd h(x_1,x_2) & = &\int
_{a_1}^{a_2} \int_{b_1}^{b_2}
h(x_1-,x_2-) \mrmd f(x_1,x_2)
\\
& &{} - \int_{a_1}^{a_2} h(x_1-,b_2)
\mrmd f_{b_2}(x_1) - \int_{b_1}^{b_2}
h(a_2,x_2-) \mrmd f_{a_2}(x_2)
\\
& &{} + \int_{a_1}^{a_2} h(x_1-,a_1)
\mrmd f_{a_1}(x_1) + \int_{b_1}^{b_2}
h(b_1,x_2-) \mrmd f_{b_1}(x_2)
\\
&&{} +f(a_2,b_2)h(a_2,b_2)-f(a_2,b_1)h(a_2,b_1)
\\
& &{} -f(a_1,b_2)h(a_1,b_2)+f(a_1,b_1)h(a_1,b_1);
\end{eqnarray*}
see Gill, van~der Laan and Wellner \cite{Gill1995}, Lemma 2.2. The
remaining part of the proof is
then similar to the proof of Lemma B.1 in Beutner, Wu and Z{\"a}hle
\cite{BeutnerWuZaehle}.
\end{pf*}

\begin{remark}\label{remarknon-implicationofboundedvariation}
(i) Again, if for $\pr$-almost every $\omega$ there are $x_{\ell
}(\omega)$ and $x_u(\omega)$ as in Remark \ref
{remarkconditionsintheonedimensioanlintbyparts}(ii), then the
conditions in part
(c) and (d) of Lemma \ref{Lemmaintegrationbypartstwodimensions}
reduce to conditions on $F$ and $g$. Furthermore, if such $x_{\ell
}(\omega)$ and $x_u(\omega)$ exist, then conditions (c) and (d) of
Lemma \ref{Lemmaintegrationbypartstwodimensions} hold whenever
$\iint|\dd g|<\infty$ 
as well as 
$\sup_{x_i\in\R}\int|\dd g_{x_i}|<\infty$ for $i=1,2$.

\mbox{}\hphantom{i}(ii) If for some weight function $\phi$ we have that $d_\phi(F_n,F)$
is $\pr$-a.s. finite, that the integral $\iint
1/(\phi(x_1)\phi(x_2)) |\dd g|(x_1,x_2)$ is finite, that
$\lim_{x_i\to\pm\infty}1/\phi(x_i)\int1/\phi(x) |\dd g_{x_i}|(x)=0$
holds for $i=1,2$, and that $g(x_1,x_2)/(\phi(x_1)\phi(x_2))$
converges to zero as $|x_1|,|x_2|\to\infty$, then again conditions
(c) and (d) of Lemma \ref{Lemmaintegrationbypartstwodimensions}
hold.

(iii) Right-continuity of $g_{x_1}$ and $g_{x_2}$, which is needed for
the integrals in part (d) of Lemma~\ref
{Lemmaintegrationbypartstwodimensions} to be well defined, is implied
by right-continuity of $g$.

\hspace*{1pt}(iv) It is worth pointing out that $g_{x_1}$ and $g_{x_2}$ being
locally of bounded variation does not imply that $g \in\BV_{\loc,\rc
}^2$; see Remark 1.1 in Beutner and Z{\"a}hle \cite{Beutnersuppl}.
Moreover, $g \in\BV
_{\loc,\rc}^2$ does not imply that $g_{x_1}$ and $g_{x_2}$ are
locally of bounded variation. Take, for example, the function $g\dvtx [0,1]
\times[0,1] \rightarrow\R$ defined by $g(0,0):=0$ and
$g(x_1,x_2):=x_1\sin(\frac{\pi}{x_1})$, $(x_1,x_2)\not=(0,0)$.

\mbox{}\hphantom{i}(v) If $g\dvtx  \R^2 \rightarrow\R$ is continuous, the partial derivative
$\partial g\slash\partial x_1$ exists and is continuous, and the mixed
partial derivative $\partial^2 g\slash(\partial x_1\,\partial x_2)$
exists and is bounded on every rectangle $\mathcal{R} \subset\R^2$,
then $g$ is locally of bounded variation; see, for instance,
Ghorpade and Limaye \cite{Inder}, Proposition 3.59.

(vi) If for all $(x_1,x_2), (y_1,y_2) \in\R^2$ with $ x_1 \leq x_2$
and $y_1 \leq y_2$ we have that
%
\begin{equation}
\label{eqremarkbimonotoneincreasing} g(x_2,y_2)+g(x_1,y_1)
\geq g(x_2,y_1)+g(x_1,y_2),
\end{equation}
then $g$ is locally of bounded variation. The same claim holds, if we
have $\leq$ instead of $\geq$ in (\ref{eqremarkbimonotoneincreasing}).
See, for instance, Ghorpade and Limaye \cite{Inder}, Proposition 1.15.
\end{remark}


\subsection{\texorpdfstring{Examples for the representation (\protect\ref{EqBZdecomposition})}
{Examples for the representation (5)}}\label{subsectionexamplesforthedecomposition}

In this section, we give some examples for set-ups under which the
representation (\ref{EqBZdecomposition}) holds. In the third,
fourth and fifth example the set-up is degenerate because there
$g_{1,F}\equiv g_{2,F}\equiv0$. Before turning to the examples, we
state two remarks including some notation needed for the examples.


\begin{remark}\label{remarkongeneratingfunctions-1}
Recall that two functions $f_1,f_2\in\BV_{\loc,\rc}^2$ generate the
same measure on $\R^2$ if
$f_1(x_1,x_2)=f_2(x_1,x_2)+h_1(x_1)+h_2(x_2)$ for some functions
$h_1,h_2\dvtx  \R\rightarrow\R$.
\end{remark}

\begin{remark}\label{remarkongeneratingfunctions-2}
For any positive measure $\mu$ on $\R$, and any measurable function
$w\dvtx \R\to\R_+$, define the measure ${\cal H}_{w,\mu}^1$ on $\R^2$ by
%
\begin{equation}
\label{remarkongeneratingfunctions-2-eq1} {\cal H}_{w,\mu}^1(A):=
\int w(x) \delta_{(x,x)}(A) \mu(\dd x),\qquad A\in{\cal B}\bigl(\R^2
\bigr),
\end{equation}
with $\delta_{(x,x)}$ the Dirac measure at $(x,x)\in\R^2$. In
particular, for ${\cal H}_{w,\mu}^1$-integrable $f\dvtx \R^2\to\R$ we
have
%
\begin{equation}
\label{defofHwmu-1} \iint f(x_1,x_2) {\cal
H}_{w,\mu}^1\bigl(\dd(x_1,x_2)\bigr) =
\int w(x)f(x,x) \mu(\dd x).
\end{equation}
The ${\cal H}_{w,\mu}^1$-measure of the area of a rectangle ${\cal R}$
intersecting the diagonal $D=\{(x,x)\dvtx x\in\R\}$ is equal to the
integral $\int_{{\cal R}_\pi}w(x)\mu(\dd x)$, where ${\cal R}_\pi$ is
the projection\vspace*{1pt} of ${\cal R}$ on one of the axes of that piece of the
diagonal $D$ that is contained in the rectangle. So one easily sees
that for every $\mathcal{R}_{a,b}=(a_1,a_2] \times(b_1,b_2]$
%
\begin{equation}
\label{defofHwmu-2} {\cal H}_{w,\mu}^1({\cal
R}_{a,b}) = \cases{\displaystyle  \int_{(\max\{a_1,b_1\},\min\{a_2,b_2\})}w(x) \mu
(\dd x), &\quad ${
\cal R}_{a,b}\cap D\not=\varnothing$,
\vspace*{2pt}\cr
0, &\quad else.}
\end{equation}

If $w\equiv1$ and $\mu$ is the Lebesgue measure $\ell$ on $\R$,
then ${\cal H}_{w,\mu}^1$ coincides with the one-dimensional Hausdorff
measure ${\cal H}^1$ in $\R^2$ restricted to the diagonal $D$ and
weighted by the constant $1/\sqrt{2}$, that is, with ${\cal H}^1(
\cdot\cap D)/\sqrt{2}$. In this case, we also write ${\cal H}_{D}^1$
instead of ${\cal H}_{\eins,\ell}^1$. As special cases of (\ref
{defofHwmu-1}) and (\ref{defofHwmu-2}) we obtain $\iint f(x_1,x_2)
{\cal H}_{D}^1(\dd(x_1,x_2)) = \int f(x,x) \mrmd x$ and
\[
{\cal H}_{D}^1({\cal R}_{a,b}) = \cases{\displaystyle  \min
\{a_2,b_2\}-\max\{a_1,b_1\}, &\quad ${
\cal R}_{a,b}\cap D\not=\varnothing$,
\vspace*{2pt}\cr
0, &\quad else.}
\]
Analogously, we let ${\cal H}_{\widetilde D}^1$ be the one-dimensional
Hausdorff measure ${\cal H}^1$ in $\R^2$ restricted to the rotated
diagonal $\widetilde D=\{(x,-x)\dvtx x\in\R\}$ and weighted by the
constant $1/\sqrt{2}$, and we note that $\iint f(x_1,x_2) {\cal
H}_{\widetilde D}^1(\dd(x_1,x_2)) = \int f(x,-x) \mrmd x$
for every ${\cal H}_{\widetilde D}^1$-integrable $f\dvtx \R^2\to\R$.
\end{remark}


\begin{example}[(Gini's mean difference)]\label{examplegini}
If $g(x_1,x_2)=|x_1-x_2|$ and $F$ has a finite first moment, then
$V_g(F)$ equals Gini's mean difference $\ex[|X_1-X_2|]$ of two
i.i.d. random variables $X_1$ and $X_2$ with df $F$. We will now
verify that, if $F$ and the estimator $F_n$ satisfy Assumptions
\ref{assumptionbasiconUFandFn} and~\ref{assumptionbasiconUFandFn-2} for
this $g$, and $d_{\phi}(F_n,F)$ is $\pr$-a.s.
finite for all $n\in\N$ and some weight function $\phi$ satisfying
$\int1/\phi(x) \mrmd x<\infty$, then the assumptions of Lemmas
\ref{Lemmaintegrationbypartsonedimension} and \ref
{Lemmaintegrationbypartstwodimensions} hold true and we have
%
\begin{equation}
\label{examplegini-eq} \dd g_{1,F}(x) = \dd g_{2,F}(x) =
\bigl(2F(x)-1\bigr) \mrmd x \quad\mbox{and}\quad \dd g(x_1,x_2) = -2{\cal
H}_D^1 \bigl(\dd(x_1,x_2) \bigr)
\end{equation}
with the notation of Remark \ref{remarkongeneratingfunctions-2}. Notice
that the left-hand side in (\ref{examplegini-eq}) shows
in particular that
the $V$-statistic corresponding to Gini's mean difference is typically
non-degenerate w.r.t. $(g,F)$ in the sense of Section \ref
{secnotionofdegeneracyandnon-dgeneracy}.

It was shown in Example 3.1 in Beutner and Z{\"a}hle \cite
{BeutnerZaehle2011} that
$g_{i,F}(x)=-\ex[X_1]+x+ 2\int_{x}^{\infty}(1-F(y)) \mrmd y$ for $i=1,2$.
Therefore, our assumption $\int1/\phi(x) \mrmd x<\infty$ implies
$\int1/\phi(x) |\dd g_{i,F}|(x)<\infty$ for $i=1,2$. From this, and
using Remark \ref{remarkconditionsintheonedimensioanlintbyparts}(iv),
we obtain the validity of assumptions (b)--(d) of
Lemma \ref{Lemmaintegrationbypartsonedimension}. It was also
established in Example 3.1 in Beutner and Z{\"a}hle \cite
{BeutnerZaehle2011} that the
left-hand side in
(\ref{examplegini-eq}) holds true. We next focus on the right-hand
side in (\ref{examplegini-eq}) and assumptions (b)--(d) of Lemma \ref
{Lemmaintegrationbypartstwodimensions}.

\textit{As for} (b): It was already shown in Example 3.1 in Beutner and
Z{\"a}hle \cite{BeutnerZaehle2011} that $g_{x_1}$ and $g_{x_2}$ are
locally of bounded variation for the above $g$. Further, for an
arbitrary rectangle $\mathcal{R}_{a,b}=(a_1,a_2] \times(b_1,b_2]$ with
$a_2 \leq b_1$ we have $ |a_2-b_2| + |a_1-b_1| - |a_1-b_2| - |a_2-b_1|
= 0$. The same holds for all rectangles $\mathcal{R}_{a,b}=(a_1,a_2]
\times(b_1,b_2]$ with $b_2 \leq a_1$. Now, consider a rectangle
$\mathcal{R}_{a,b}=(a_1,a_2] \times(b_1,b_2]$ with $a_2 \geq b_2 > a_1
\geq b_1$. Then $|a_2-b_2| + |a_1-b_1| - |a_1-b_2| - |a_2-b_1| =
2(a_1-b_2) < 0$. Similar inequalities hold for the remaining cases,
that is, for $a_2 \geq b_2 > b_1 > a_1$, $b_2 \geq a_2 > b_1 \geq a_1$, and
$b_2 \geq a_2 > a_1 > b_1$. Hence, Remark \ref
{remarknon-implicationofboundedvariation}(vi) implies $g
\in\BV_{\loc,\rc}^2$.

\textit{As for} (c): Let $\mu_g$ denote the measure generated by Gini's
mean difference kernel. We just have seen that
$\mu_g(\mathcal{R}_{a,b})=0$ for an arbitrary half-open rectangle
$\mathcal{R}_{a,b}$ not intersecting the diagonal. Moreover, we have
seen that $\mu_g(\mathcal{R}_{a,b})=2(a_1-b_2)$ when $b_1 \leq a_1 <
b_2 \leq a_2$. Taking the other possibilities mentioned above into
account, we find that
\[
\mu_g(\mathcal{R}_{a,b}) = \cases{\displaystyle  2\bigl(\max
\{a_1,b_1\}-\min\{a_2,b_2\}\bigr),
&\quad ${\cal R}_{a,b}\cap D\not=\varnothing$,
\vspace*{2pt}\cr
0, &\quad else.}
\]
Thus, in view of Remarks \ref{remarkongeneratingfunctions-1}--\ref
{remarkongeneratingfunctions-2}, the measure $\mu_g$
generated by Gini's mean difference kernel differs from ${\cal
H}_{D}^1$ only by the sign and the factor $2$, that is, the right-hand
side in (\ref{examplegini-eq}) holds. In view of Remark \ref
{remarknon-implicationofboundedvariation}(ii), this implies
condition (c) of Lemma \ref{Lemmaintegrationbypartstwodimensions},
because we assumed $d_{\phi}(F_n,F)<\infty$ $\pr$-a.s.
for all $n\in\N$ and some weight function $\phi$ with $\int1/\phi
(x) \mrmd x<\infty$.

\textit{As for} (d): It was shown in Example 3.1 in
Beutner and Z{\"a}hle \cite{BeutnerZaehle2011} that $\dd g_{x_i}^+(x) =
\eins_{(x_i,\infty]}(x) \mrmd x$ and $\dd g_{x_i}^-(x) =
\eins_{[-\infty,x_i]}(x) \mrmd x$ for $i=1,2$. From this and the obvious
convergence of $|x_1-x_2|/(\phi(x_1)\phi(x_2))$ to zero as
$|x_1|,|x_2|\to\infty$, along with Remark \ref
{remarknon-implicationofboundedvariation}(ii), it can be deduced
easily that all limits in condition (d) of Lemma \ref
{Lemmaintegrationbypartstwodimensions} exist and equal zero $\pr$-a.s.
\end{example}


\begin{example}[(Variance)]\label{examplevariance}
If $g(x_1,x_2)=\frac{1}{2}(x_1-x_2)^2$ and $F$ has a finite second
moment, then $V_g(F)$ equals the variance of $F$. We will now
verify that, if $F$ and the estimator $F_n$ satisfy Assumptions
\ref{assumptionbasiconUFandFn} and \ref{assumptionbasiconUFandFn-2} for
this $g$, and $d_{\phi}(F_n,F)$ is $\pr$-a.s.
finite for all $n\in\N$ and some weight function $\phi$ satisfying
$\int|x|/\phi(x) \mrmd x<\infty$, then the assumptions of Lemmas
\ref{Lemmaintegrationbypartsonedimension} and \ref
{Lemmaintegrationbypartstwodimensions} hold true and we have
%
\begin{equation}
\label{examplevariance-eq-two} \dd g_{1,F}(x) = \dd g_{2,F}(x) =
\bigl(x-\ex[X_1]\bigr) \mrmd x \quad\mbox{and}\quad \dd g(x_1,x_2)
= - \dd x_1 \mrmd x_2.
\end{equation}
Notice that the left-hand side of (\ref{examplevariance-eq-two})
shows in particular that
the $V$-statistic corresponding to the variance is typically
non-degenerate w.r.t. $(g,F)$ in the sense of Section \ref
{secnotionofdegeneracyandnon-dgeneracy}.

It was already verified in Example 3.2 in Beutner and Z{\"a}hle \cite
{BeutnerZaehle2011}
that $g_{i,F}(x)=\frac{1}{2}x^2- x \ex[X_1]+\frac{1}{2}\ex[X_1^2]$
for $i=1,2$. Therefore, our assumption $\int|x|/\phi(x) \mrmd x<\infty$
implies $\int1/\phi(x) |\dd g_{i,F}|(x)<\infty$ for $i=1,2$. Thus,
$g_{i,F} \in\BV_{\loc,\rc}$ and condition (c) of Lemma \ref
{Lemmaintegrationbypartsonedimension} follows at once. Moreover,
assumption (d) of Lemma \ref{Lemmaintegrationbypartsonedimension} holds
by Remark \ref{remarkconditionsintheonedimensioanlintbyparts}(v). It
was also established in Example
3.2 in Beutner and Z{\"a}hle \cite{BeutnerZaehle2011} that the
left-hand side of (\ref
{examplevariance-eq-two})
holds true. We next focus on the right-hand side of (\ref
{examplevariance-eq-two}) and
assumptions (b)--(d) of Lemma \ref{Lemmaintegrationbypartstwodimensions}.

\textit{As for} (b): It was already shown in Example 3.2 in
Beutner and Z{\"a}hle \cite{BeutnerZaehle2011} that $g_{x_1}$ and
$g_{x_2}$ are locally
of bounded variation for the above $g$. Moreover, we have $g \in
\BV_{\loc,\rc}^2$ by Remark \ref
{remarknon-implicationofboundedvariation}(v).

\textit{As for} (c) \textit{and} (d): Notice that $g(x_1,x_2)= -x_1x_2+
\frac{1}{2}x_1^2+ \frac{1}{2}x_2^2$ and recall Remark \ref
{remarkongeneratingfunctions-1}. Thus, up to the sign, the measure
generated by the variance kernel is equal to the Lebesgue measure
on~$\R^2$, that is, the right-hand side of
(\ref{examplevariance-eq-two}) holds. Moreover, it was verified in
Example 3.2 in Beutner and Z{\"a}hle \cite {BeutnerZaehle2011} that
$\dd g_{x_i}^+(x)= (x-x_i) \eins_{(x_i,\infty]}(x) \mrmd x$ and $\dd
g_{x_i}^-(x) = (x_i-x) \eins_{[-\infty,x_i]}(x) \mrmd x$ for $i=1,2$.
Thus, we see from Remark
\ref{remarknon-implicationofboundedvariation}(ii), that conditions (c)
and (d) of Lemma \ref {Lemmaintegrationbypartstwodimensions} hold.
\end{example}


Now, let us turn to some examples where the linear part of the
representation (\ref{EqBZdecomposition}) vanishes.

\begin{example}[(Gini's mean difference, degenerate case)]\label{exampleginidegenrate}
The $V$-statistic $V_g(\hat F_n)$ corresponding to Gini's mean difference
kernel $g(x_1,x_2)=|x_1-x_2|$ is degenerate w.r.t. $(g,F)$ for any df
$F$ that assigns probability $1/2$ to two points in $\R$. Indeed, we
know from (\ref{examplegini-eq}) that\vadjust{\goodbreak}
$\dd g_{1,F}(x)=\dd g_{2,F}(x)=(2F(x)-1)\mrmd x$ for all $x \in\R$, so that
$\int(\hat{F}_n(x)-F(x)) \mrmd g_{i,F}(x)=0$, $i=1,2$, in this case.
Recall that $\hat{F}_n$ refers to the empirical df, and notice that
the assumptions of Lemma \ref{Lemmaintegrationbypartstwodimensions}
trivially hold, because $(\hat{F}_n-F)(x)$ equals zero for
$x$ small and large enough, respectively.
\end{example}


\begin{example}[(Goodness-of-fit test)]\label{gof}
For a given df $F_0$ and any measurable (weight) function $w\dvtx \R\to\R
_+$, the weighted Cram\'er--von Mises test statistic
$T_n^0:=\int w(x)(\hat F_n(x)-F_0(x))^2 \mrmd F_0(x)$
was introduced for testing the null hypothesis $F=F_0$, and coincides
with the classical Cram\'er--von Mises test statistic for $w\equiv1$.
The test statistic $T_n^0$ can be expressed as $V$-statistic $V_g(\hat
F_n)$ with kernel
%
\begin{equation}
\label{tfs-eq-1} g(x_1,x_2):= \int w(x) \bigl(
\eins_{[x_1,\infty)}(x)-F_0(x) \bigr) \bigl(\eins_{[x_2,\infty)}(x)-F_0(x)
\bigr) \mrmd F_0(x).
\end{equation}
We will now verify that, if $F_0$ is continuous and satisfies
Assumptions \ref{assumptionbasiconUFandFn} and \ref
{assumptionbasiconUFandFn-2} for this $g$ and if the integral $\int
w(x)
\mrmd F_0(x)$ is finite, then under the null hypothesis $F=F_0$ the
assumptions of Lemmas \ref{Lemmaintegrationbypartsonedimension}
and \ref{Lemmaintegrationbypartstwodimensions} hold true and we have
%
\begin{equation}
\label{tfs-eq-2} g_{1,F_0} \equiv g_{2,F_0} \equiv0
\quad\mbox{and}\quad
\dd g(x_1,x_2) = {\cal H}_{w,\dd F_0}^1
\bigl(\dd(x_1,x_2) \bigr)
\end{equation}
with the notation of Remark \ref{remarkongeneratingfunctions-2}.

The left-hand side in (\ref{tfs-eq-2}) follows by using Fubini's theorem.
Hence, the assumptions of Lemma \ref
{Lemmaintegrationbypartsonedimension} trivially hold in this case. We
next focus on the right-hand
side in (\ref{tfs-eq-2}) and assumptions (b)--(d) of Lemma \ref
{Lemmaintegrationbypartstwodimensions}. From (\ref{tfs-eq-1}), we
easily see that under our assumptions the sections $g_{x_1}$ and
$g_{x_2}$ are right-continuous and locally of bounded variation.
Moreover, the right-hand side in (\ref{tfs-eq-2}) is known from
Dehling and Taqqu \cite{DehlingTaqqu1991}, Example 3, and so it is
apparent that $g \in
\BV_{\loc,\rc}^2$. Hence, (b) holds. From (\ref{tfs-eq-2}),
(\ref{defofHwmu-1}), and the fact that we assumed $\int
w(x) \mrmd F_0(x)$ to exist, we immediately obtain (c). Further, we note
that $\dd g_{x_i}^+(x) = w(x)F_0(x) \mrmd F_0(x)$ and $\dd g_{x_i}^-(x) =
w(x)\eins_{[x_i,\infty)}(x) \mrmd F_0(x)$ for $i=1,2$. From this and
Remark \ref{remarknon-implicationofboundedvariation}(i), which
can be applied since $\hat{F}_n$ is the empirical df, it can now be
easily deduced that all limits in condition (d) of Lemma \ref
{Lemmaintegrationbypartstwodimensions} exist and equal zero $\pr$-a.s.
\end{example}


\begin{example}[(Test for symmetry)]\label{tfs}
The statistic $T_n:=\int_0^\infty(\hat F_n(-t)-[1-\hat F_n(t-)])^2 \mrmd t$
is often used for testing symmetry of $F$ about zero; cf.
Arcones and Gin{\'e} \cite{Arcones92}, Example 5.1. Using Fubini's
theorem, more
precisely Theorem 1.15 in Mattila \cite{Mattila1995} with $d\mu=d\hat
F_n\times d\hat F_n$ and its analogue for negative integrands, $T_n$
can be expressed as $V$-statistic $V_g(\hat F_n)$ with kernel
%
\begin{equation}
\label{tfsym-eq-1} g(x_1,x_2):= \bigl(|x_1|
\wedge|x_2|\bigr) (\eins_{\{x_1,x_2>0\}
}+\eins_{\{x_1,x_2<0\}}-
\eins_{\{x_1>0, x_2<0\}}-\eins_{\{x_1<0,
x_2>0\}} ).
\end{equation}
We will now verify that, if $F$ satisfies Assumptions
\ref{assumptionbasiconUFandFn} and \ref{assumptionbasiconUFandFn-2} for
this $g$, and $d_{\phi}(\hat F_n,F_0)$ is $\pr
$-a.s. finite for
all $n\in\N$ and some symmetric weight function $\phi$ with $\int
1/\phi(x) \mrmd x<\infty$, then under the null hypothesis that $F$ be
symmetric about zero the assumptions of Lemmas \ref
{Lemmaintegrationbypartsonedimension} and \ref
{Lemmaintegrationbypartstwodimensions} hold true and we have
%
\begin{equation}
\label{tfsym-eq-2-two} g_{1,F} \equiv g_{2,F} \equiv0
\quad\mbox{and}\quad \dd g(x_1,x_2) = {\cal H}_{D}^1
\bigl(d(x_1,x_2) \bigr)-{\cal H}_{\widetilde D}^1
\bigl(d(x_1,x_2) \bigr)
\end{equation}
with the notation of Remark
\ref{remarkongeneratingfunctions-2}.\vadjust{\goodbreak}

The left-hand side of (\ref{tfsym-eq-2-two}) is obvious (under
the null
hypothesis), and so the assumptions of Lemma \ref
{Lemmaintegrationbypartsonedimension} trivially hold in this case. We
next focus
on the right-hand side of (\ref{tfsym-eq-2-two}) and assumptions
(b)--(d) of Lemma
\ref{Lemmaintegrationbypartstwodimensions}. From (\ref{tfsym-eq-1}) we
easily see that the sections $g_{x_1}$ and $g_{x_2}$ are
locally of bounded variation. Moreover, we have $\dd g^+(x_1,x_2)={\cal
H}_{D}^1(\dd(x_1,x_2))$ and $\dd g^-(x_1,x_2)={\cal H}_{\widetilde
D}^1(\dd(x_1,x_2))$, and so it is apparent that $g \in
\BV_{\loc,\rc}^2$, and that\vspace*{1pt} the right-hand side of (\ref
{tfsym-eq-2-two}) holds.
Hence, (b) holds. From the right-hand side of (\ref{tfsym-eq-2-two}),
(\ref{defofHwmu-1}) and our assumptions we immediately obtain (c). Further,
it can be checked easily that for $i=1,2$
\[
\dd g_{x_i}(x) = \cases{\displaystyle  \eins_{[0,x_i]}(x) \mrmd x-
\eins_{[-x_i,0]}(x) \mrmd x, &\quad $x_i>0$,
\vspace*{2pt}\cr
\displaystyle -\eins_{[x_i,0]}(x)
\mrmd x+\eins_{[0,-x_i]}(x) \mrmd x, &\quad $x_i<0$.}
\]
From this, our assumptions and 
the fact that $\hat{F}_n$ is the empirical df, it can now be easily deduced
that the first two limits in condition (d) of Lemma
\ref{Lemmaintegrationbypartstwodimensions} exist and equal zero
$\pr$-a.s. Using our assumption $\int1/\phi(x) \mrmd x<\infty$ and
Remark \ref{remarknon-implicationofboundedvariation}(ii), it
follows that the third limit in condition (d) of Lemma \ref
{Lemmaintegrationbypartstwodimensions} exists and equals zero
$\pr$-a.s.
\end{example}


\subsection{Weak (central) limit theorems}\label{secweaklimittheorems}

In this section, we give a tool for deriving the asymptotic
distribution of $V$-statistics, which is suitable for independent and
weakly dependent data. Moreover, in some particular cases it also
yields a non-trivial asymptotic distribution for strongly dependent
data. Recall that $(\V,d)$ is some metric subspace of $\D$, and that
${\cal V}={\cal D}\cap\V$.

\begin{theorem}\label{theoremcontmappingforu-stat}
Let $(a_n)$ be a sequence in $(0,\infty)$ with $a_n\to\infty$.
Assume that:
\begin{itemize}[(a)]
\item[(a)] the assumptions of Lemmas \ref
{Lemmaintegrationbypartsonedimension} and \ref
{Lemmaintegrationbypartstwodimensions} are
fulfilled,
\item[(b)] on $\V$ the functions $\Phi_{i,g}$, $i=1,2,3$, defined in
(\ref{eqdefinitionofPhii}), are well-defined, $({\cal V},{\cal
B}(\R))$-measurable and $d$-continuous,
\item[(c)] the process $a_n(F_n-F)$ is a random element of $(\V,{\cal
V})$ for all $n\in\N$, and there is some random element $B^{\circ}$ of
$(\V,{\cal V})$ such that $\pr[B^{\circ}\in S]=1$ for some
$d$-separable $S\in{\cal V}$ consisting of $({\cal V},d)$-completely
regular points (in the sense of Definition \textup{IV.2.6} in Pollard \cite
{Pollard1984}) only, and
%
\begin{equation}
\label{conditiononBcirc} a_n(F_n-F)\stackrel{\mathsf d} {
\longrightarrow}B^{\circ} \qquad\mbox{in $(\V,{\cal V},d)$}.
\end{equation}
\end{itemize}
Then the following assertions hold:
\begin{itemize}[(ii)]
\item[(i)] We always have
%
\begin{equation}
\label{theoremcontmappingforu-stat-eq1} a_n\bigl(V_g(F_n)-V_g(F)
\bigr) \stackrel{\mathsf d} {\longrightarrow} -\sum_{i=1}^2
\int B^\circ(x-) \mrmd g_{i,F}(x) \qquad\mbox{in $\bigl(\R,{\cal B}(\R)
\bigr)$}.
\end{equation}
\item[(ii)] If the $V$-statistic $V_g(F_n)$ is degenerate w.r.t.
$(g,F)$, then we additionally have
%
\begin{equation}
\label{theoremcontmappingforu-stat-eq2} a_n^2
\bigl(V_g(F_n)-V_g(F)\bigr) \stackrel{\mathsf d}
{\longrightarrow} \iint B^{\circ}(x_1-)B^{\circ}(x_2-)
\mrmd g(x_1,x_2) \qquad\mbox{in $\bigl(\R,{\cal B}(\R)\bigr)$}.
\end{equation}
\end{itemize}
\end{theorem}

\begin{pf}
(i) Under condition (a), we can apply Lemmas \ref
{Lemmaintegrationbypartsonedimension} and \ref
{Lemmaintegrationbypartstwodimensions} to obtain the
representation~(\ref{EqBZdecompositionwithPhis}). From the Continuous
Mapping theorem, which is applicable by conditions (b) and (c), we
obtain that $\Phi _{i,g}(a_n(F_n-F)) \stackrel{\mathsf
d}{\longrightarrow} -\int B^\circ (x-) \mrmd g_{i,F}(x), i=1,2$. From
Slutsky's lemma, we obtain that $\sqrt{a_n}(F_n-F)$ converges to zero
in probability, and so, according to the Continuous Mapping theorem,
$\Phi_{3,g}(\sqrt{a_n}(F_n-F))$ converges in probability to zero as
well. Applying once again Slutsky's lemma finishes the proof of part~(i).

(ii) If the $V$-statistic is degenerate, then we obtain analogously by
applying Lemmas \ref{Lemmaintegrationbypartsonedimension}
and~\ref{Lemmaintegrationbypartstwodimensions} that
$a_n^2(V_g(F_n)-V_g(F)) = \Phi_{3,g}(a_n(F_n-F))$. The result is now
immediate from the Continuous Mapping theorem.
\end{pf}

\begin{remark}\label{remarkoncontinuityofVg}
If, for some weight function $\phi$, the integral $\int
1/\phi(x) |\dd g_{i,F}|(x)$ is finite for $i=1,2$, then the
mappings $\Phi_{1,g}$ and $\Phi_{2,g}$ are obviously
$d_\phi$-continuous. Moreover, if the integral
$\int1/(\phi(x_1)\phi(x_2)) |\dd g|(x_1,x_2)$ is finite, then
$\Phi_{3,g}$ is $d_{\phi}$-continuous, too.
\end{remark}

\begin{example}[(I.i.d. data)]\label{exampleiiddata}
Suppose $X_1,X_2,\ldots$ are i.i.d. with df $F$, and let $\phi$ be a
weight function. If\vspace*{-1pt} $\int\phi^2\mrmd F<\infty$, then Theorem 6.2.1 in
Shorack and Wellner \cite{ShorackWellner1986} shows that
$\sqrt{n}(\hat F_n-F)\stackrel{\mathsf d}{\longrightarrow} B_F^\circ$
(in $(\D_\phi,{\cal D}_\phi,d_\phi)$) with $\hat F_n$ the empirical
df of $X_1,\ldots,X_n$ and $B_F^\circ$ an $F$-Brownian bridge, that
is, a centered Gaussian process with covariance function $\Gamma
(x,y)=F(x\wedge y)\overline{F}(x\vee y)$.
\end{example}

\begin{example}[(Weakly dependent data)]\label
{exampleweaklydependentdata}
Let $(X_i)$ be $\alpha$-mixing with mixing coefficients satisfying
$\alpha(n)=\RMO(n^{-\theta})$ for some $\theta>1+\sqrt{2}$, and
let $\lambda\ge0$. If\vspace*{1pt} $F$ has a finite $\gamma$-moment for some
$\gamma>(2\theta\lambda)/(\theta-1)$, then\vspace*{-1pt} it can easily be deduced
from Theorem 2.2 in Shao and Yu \cite{ShaoYu1996} that $\sqrt{n}(\hat
F_{n}-F)\stackrel{\mathsf d}{\longrightarrow}\widetilde B_F^\circ$ (in
$(\D_{\phi_\lambda},{\cal D}_{\phi_\lambda},d_{\phi_\lambda})$) with
$\widetilde B_F^\circ$ a continuous centered Gaussian process with
covariance function
$\Gamma(y_0,y_1)=F_0(y_0\wedge y_1)\overline{F}_0(y_0\vee
y_1)+\sum_{i=0}^1\sum_{k=2}^{\infty}\covi(\eins_{\{X_1 \le y_i\}},
\eins_{\{X_k \le y_{1-i}\}})$; see Section 3.3 in
Beutner and Z{\"a}hle \cite{BeutnerZaehle2010a}. If $(X_i)$ is even
$\beta$- or
$\rho$-mixing, then Lemma 4.1 in Chen and Fan \cite{ChenFan2006} and
Theorem 2.3 in
Shao and Yu \cite{ShaoYu1996} ensure that the mixing condition can be
relaxed; see
also Section 3.2 in Beutner and Z{\"a}hle \cite{BeutnerZaehle2011}.
\end{example}

Part (ii) of Theorem \ref{theoremcontmappingforu-stat} also leads to
some interesting corollaries that provide
convenient alternatives to the common approach to derive the
asymptotic distribution of degenerate $U$-statistics. Before
presenting them, it is worth recalling that the common approach to
derive the asymptotic distribution of degenerate $U$-statistics is
based on a series expansion of the kernel $g$ of the form $g(x_1,x_2)
= \sum_{k=1}^{\infty} \lambda_k \psi_k(x_1)\psi_k(x_2)$,
where the $\lambda_k$ are real numbers and the $\psi_k$ are an
orthonormal sequence; see, for example, Serfling \cite{Serfling1980}, Section
5.5. The $\lambda_k$ and the $\psi_k$ are the eigenvalues and
eigenfunctions, respectively, of the operator $A\dvtx  L^2(\R,{\cal B}(\R
),\pr) \rightarrow L^2(\R,{\cal B}(\R),\pr)$ defined by
$A(h(x_1)) = \int g(x_1,x_2)h(x_2) \mrmd F(x_2)$.
The eigenvalues arise in the asymptotic distribution of
$n(U_{g,n}-V_g(F))$ which, in the i.i.d. case, is given by
$\sum_{i=1}^{\infty} \lambda_i (\xi_i^2-1)$, where the $\xi_i$ are
independent and have a $\chi^2$-distribution with $1$ degree of
freedom.

\begin{corollary}\label{corollarydegenerateproductmeasure}
Assume that the conditions of Theorem \ref{theoremcontmappingforu-stat}
hold and that we are in case \textup{(ii)} of this theorem. Moreover,
let the (possibly signed) measure generated by $g$ on $\R^2$ be equal
to the product measure of (possible signed) measures $\nu_1$ and $\nu
_2$. Then
\[
a_n^2\bigl(V_g(F_n)-V_g(F)
\bigr)\stackrel{\mathsf d} {\longrightarrow} \biggl(\int B^{\circ}(x-)\mrmd
\nu_1(x) \biggr) \biggl(\int B^{\circ}(x-)\mrmd \nu_2(x)
\biggr) \qquad\mbox{in $\bigl(\R,{\cal B}(\R)\bigr)$}.
\]
\end{corollary}

The next example shows that two well known kernels are covered by
Corollary \ref{corollarydegenerateproductmeasure}.

\begin{example}\label{exampledegneratevariaceandsquared}
(i) The variance kernel $g(x_1,x_2)=\frac{1}{2}(x_1-x_2)^2= -x_1x_2+
\frac{1}{2}x_1^2+ \frac{1}{2}x_2^2$ is degenerate if and only if the
fourth central moment equals the squared second central moment; see,
for example, van der Vaart \cite{vanderVaart1998}, Example 12.12.
Moreover, in
Example \ref{examplevariance} we have seen that the measure generated
by the variance kernel coincides with the negative Lebesgue measure on
$\R^2$. So, in the degenerate case, the variance kernel can be treated
by means of Corollary \ref{corollarydegenerateproductmeasure}.

(ii) The kernel $g(x_1,x_2)=x_1x_2$, which corresponds to the
characteristic $\ex[X_1]^2$ and which is degenerate if the first
moment equals zero, obviously generates the Lebesgue measure. In
particular, up to the sign, it generates the same measure on $\R^2$ as
the variance kernel. So, in the degenerate case, this kernel can be
treated by means of Corollary \ref{corollarydegenerateproductmeasure}
as well. Of course, for the corresponding $V$-statistics the
asymptotic distributions can be derived differently, but the continuous
mapping approach reveals an interesting relation to the variance kernel.
\end{example}


\begin{corollary}\label{corollarydegeneratedirackernel}
Assume that the conditions of Theorem \ref{theoremcontmappingforu-stat}
hold and that we are in case \textup{(ii)} of this theorem. Moreover,
let the measure generated by $g$ be given by ${\cal H}_{w,\mu}^1$
defined in (\ref{remarkongeneratingfunctions-2-eq1}). Then
\begin{eqnarray*}
a_n^2\bigl(V_g(F_n)-V_g(F)
\bigr) \stackrel{\mathsf d} {\longrightarrow} \int w(x) \bigl(B^{\circ}(x-)
\bigr)^2 \mu(\dd x) \qquad\mbox{in $\bigl(\R,{\cal B}(\R)\bigr)$}.
\end{eqnarray*}
\end{corollary}

Here are two examples to which Corollary \ref
{corollarydegeneratedirackernel} can be applied.

\begin{example}\label{exampleGiniandCramer-vonMisesdegenerate}
(i) In Example \ref{exampleginidegenrate}, we have seen that Gini's
mean difference is degenerate for a df that assigns probability $1
\slash2$ to two points in $\R$. Further, from Example \ref{examplegini}
we also know that the measure generated by $g$ differs from
${\cal H}_{\eins,\ell}^1={\cal H}_{D}^1$ only by a constant factor.
So, Corollary~\ref{corollarydegeneratedirackernel} can be applied.

(ii) In Example \ref{gof}, we have seen that the measure generated by
the kernel $g$ of the Cram\'{e}r--von Mises statistic (cf. (\ref
{tfs-eq-1})) equals the measure ${\cal H}_{w,\dd F}^1$. Thus, the
Cram\'{e}r--von Mises statistic can also be treated by means of
Corollray \ref{corollarydegeneratedirackernel}. For the particular
case $w \equiv1$ see also van der Vaart \cite{vanderVaart1998},
Corollary 19.21.
\end{example}

Although, the asymptotic distributions in Example \ref
{exampleGiniandCramer-vonMisesdegenerate} can be derived differently,
the two
examples given are appealing from a structural point of view.


\subsection{Strong limit theorems}\label{secstronglimittheorems}

In this section, we focus on almost sure convergence of the plug-in
estimator $V_g(F_n)$ to $V_g(F)$. Assume that the representation (\ref
{EqBZdecomposition}) holds, that the mapping $\Phi\dvtx \V\to\R$,
$f\mapsto(\sum_{i=1}^2\Phi_{i,g}(f)+\Phi_{3,g}(f))$, is
$d$-continuous at the null function, and that $F_n-F$ converges $\pr
$-a.s. to the null function w.r.t. $d$. Then we immediately obtain
almost sure convergence of $V_g(F_n)$ to $V_g(F)$. From the following
obvious theorem, we can even deduce the rate of convergence. By local
$\beta$-H\"{o}lder $d$-continuity of a functional $\Phi\dvtx \V\to\R$
at $f$, we mean that $|\Phi(f_n)-\Phi(f)|=\RMO(d(f_n,f)^\beta)$
for each sequence $(f_n)\subset\V$ with $d(f_n,f)\to0$.

\begin{theorem}\label{theoremmzslln}
Let $\phi\dvtx \R\to[1,\infty)$ be some weight function, let $d$ be
homogeneous, and assume that:
\begin{itemize}[(a)]
\item[(a)] the assumptions of Lemmas \ref
{Lemmaintegrationbypartsonedimension} and \ref
{Lemmaintegrationbypartstwodimensions} are
fulfilled,
\item[(b)] on $\V$ the functions $\Phi_{i,g}$, $i=1,2,3$, defined in
(\ref{eqdefinitionofPhii}), are well defined and $({\cal V},{\cal
B}(\R))$-measurable, and the function $\sum_{i=1}^3\Phi_{i,g}$ is
locally $\beta$-H\"older $d$-continuous at the null function for some
$\beta>0$, 
\item[(c)] the process $F_n-F$ is a random element of $(\V,{\cal V})$
for all $n\in\N$, and, for some sequence $(a_n)$ in $(0,\infty)$,
%
\begin{equation}
\label{theoremmzslln-eq} a_n d(F_n-F,0)\longrightarrow0\qquad
\pr\mbox{-a.s.}
\end{equation}
\end{itemize}
Then
\[
a_n^\beta\bigl(V_g(F_n)-V_g(F)
\bigr)\longrightarrow0\qquad \pr\mbox{-a.s.}
\]
\end{theorem}

\begin{remark}
If, for some weight function $\phi$, the integral
$\int\phi(x) |\dd g_{i,F}|(x)$ is finite for $i=1,2$, then the
functionals $\Phi_{1,g}$ and $\Phi_{2,g}$ are obviously locally
$1$-H\"older $d_\phi$-continuous at the null function. Moreover, if
the integral $\int1/(\phi(x_1)\phi(x_2)) |\dd g|(x_1,x_2)$ is finite,
then the functional $\Phi_{3,g}$ is obviously $2$-H\"older
$d_{\phi}$-continuous at the null function. Thus, in this case the
functional $\sum_{i=1}^3\Phi_{i,g}$ is locally $1$-H\"older
$d_{\phi}$-continuous at the null function, and the rate of
convergence of degenerate $V$-statistics w.r.t. $(g,F)$, that is, of
$V$-statistics with $\sum_{i=1}^2\Phi_i(F_n-F)=0$, is twice the rate
of non-degenerate $V$-statistics w.r.t. $(g,F)$.
\end{remark}



\begin{example}\label{generalizedGC-independent-mz}
(i) Let $\phi$ be any weight function, and $r\in[0,\frac{1}{2})$. If
the sequence $(X_i)$ is i.i.d. and $\int\phi^{1/(1-r)}\mrmd F<\infty$, then
(\ref{theoremmzslln-eq}) hold for the weighted sup-metric
$d:=d_\phi$, the empirical df $F_n:=\hat F_n$, and $a_n:=n^r$; cf.
Andersen, Gin{\'e} and Zinn \cite{AndersenGineZinn1988}, Theorem 7.3.

\mbox{}\hphantom{i}(ii) Suppose that $\int\phi \mrmd F<\infty$. Further suppose that $(X_i)$
is $\alpha$-mixing with mixing coefficients $\alpha(n)$, let
$\alpha(t):=\alpha(\lfloor t\rfloor)$ be the c\`adl\`ag extension of
$\alpha(\cdot)$ from $\N$ to $\R_+$, and assume that $\int_0^1
\log(1+\alpha^\rightarrow(s/2) ) \overline{G}
{}^\rightarrow(s) \mrmd s<\infty$
for $\overline G:=1-G$, where $G$ denotes the df of $\phi(X_1)$ and
$\overline{G} {}^\rightarrow$ the right-continuous inverse of
$\overline{G}$. It was shown in Z{\"a}hle \cite{Zaehle2012} that, under the
imposed assumptions, (\ref{theoremmzslln-eq})~holds for the
weighted sup-metric $d:=d_\phi$, the empirical df $F_n:=\hat F_n$, and
$a_n:=1$. Notice that the integrability condition above holds in
particular if $\ex[\phi(X_1)\log^+\phi(X_1)]<\infty$ and
$\alpha(n)=\RMO(n^{-\vartheta})$ for some arbitrarily small
$\vartheta>0$; cf. Rio \cite{Rio1994}, Application 5, page 924.

(iii) Suppose that the sequence $(X_i)$ is $\alpha$-mixing with mixing
coefficients $\alpha(n)$. Let $r\in[0,\frac{1}{2})$ and assume that
$\alpha(n)\le K n^{-\vartheta}$ for all $n\in\N$ and some constants
$K>0$ and $\vartheta>2r$. Then (\ref{theoremmzslln-eq}) holds for
the uniform sup-metric $d:=d_\infty$, the empirical df $F_n:=\hat F_n$,
and $a_n:=n^r$; cf. Z{\"a}hle~\cite{Zaehle2012}.
\end{example}


\section{\texorpdfstring{The use of the representation (\protect\ref{EqBZdecomposition}) for linear long-memory sequences}
{The use of the representation (5) for linear long-memory sequences}}\label{sectionapplication}

As indicated in the \hyperref[sectionintroduction]{Introduction} and in Section \ref
{secnotionofdegeneracyandnon-dgeneracy}, for sequences exhibiting long-range
dependence it may happen that the linear part of the von Mises
decomposition degenerates only asymptotically. In such a case,
Theorem \ref{theoremcontmappingforu-stat} may not yield a
non-central limit theorem; for more details see the discussion below
just after the proof of Theorem \ref{nclt}. Nevertheless,
representation (\ref{EqBZdecomposition}), the Continuous Mapping
theorem, and an ``expansion'' of the empirical process will lead to
a general result to derive non-central limit theorems for $U$- and
$V$-statistics based on linear long-memory sequences. Thus, in this
section, we shall consider a linear process exhibiting long-range
dependence (strong dependence), that is,
\[
X_t:= \sum_{s=0}^\infty
a_s \varepsilon_{t-s},\qquad t\in\N,
\]
%
where $(\varepsilon_i)_{i\in\Z}$ are i.i.d. random variables on
some probability space $(\Omega,{\cal F,\pr})$ with zero mean and
finite variance, and the coefficients $a_s$ satisfy $\sum_{s=0}^\infty
a_s^2<\infty$ (so that $(X_t)_{t\in\N}$ is an $L^2$-process) and
decay sufficiently slowly so that $\sum_{t=1}^\infty|{\covi}
(X_1,X_t)|=\infty$. The latter divergence gives the precise meaning to
the attribute long-range dependence. Notice that if $\varepsilon_1$
has a finite $p$th moment for some $p\ge2$, then the same holds for
$X_1$. As before, we denote by $F$ the df of the $X_t$.

For $n\in\N$ and $p\in\N_0$, assume that the $p$th moment of $F$ is
finite and that $F$ can be differentiated at least $p$ times. Denote
the $j$th derivative of $F$ by $F^{(j)}$, $j=0,\ldots,p$, with the
convention $F^{(0)}=F$, and define a stochastic process ${\cal
E}_{n,p;F}$ with index set $\R$ by
%
\begin{eqnarray}
\label{EqVarepsprocess} {\cal E}_{n,p;F}(\cdot) &:=& \hat
F_n(\cdot) - \sum_{j=0}^p(-1)^{j}
F^{(j)}(\cdot) \Biggl(\frac{1}{n}\sum
_{i=1}^nA_{j;F}(X_i) \Biggr)
\nonumber\\[-8pt]\\[-8pt]
& = & \hat F_n(\cdot)-F(\cdot) - \sum_{j=1}^p(-1)^{j}
F^{(j)}(\cdot) \Biggl(\frac{1}{n}\sum
_{i=1}^nA_{j;F}(X_i) \Biggr),\nonumber
\end{eqnarray}
where $A_{j;F}$ denotes the $j$th order Appell polynomial associated
with $F$, and we use the convention $\sum_{j=1}^0(\cdots):=0$.
Recall that these Appell polynomials are defined by $A_{0;F}(x):=1$ and
for $j=1,\ldots,p$ recursively by the characteristic conditions
\[
\frac{\dd}{\dd x}A_{j;F}(x) = j A_{j-1;F}(x) \quad\mbox{and}\quad \int
A_{j;F}(y) \mrmd F(y) = 0.
\]
In particular,
${\cal E}_{n,0;F}(\cdot)=(\hat F_n(\cdot)-F(\cdot))$ and ${\cal
E}_{n,1;F}(\cdot)=(\hat F_n(\cdot)-F(\cdot))+F^{(1)}(\cdot)(\frac
{1}{n}\sum_{i=1}^n X_i)$.
For $p\in\N$, we obviously have
%
\begin{equation}
\label{generalizedempproc} {\cal E}_{n,p-1;F}(\cdot) = {\cal
E}_{n,p;F}(\cdot) + (-1)^{p} F^{(p)}(\cdot) \Biggl(
\frac{1}{n}\sum_{i=1}^n
A_{p;F}(X_i) \Biggr),
\end{equation}
and we note that under a suitable re-scaling the limit in distribution,
$Z_{p,\beta}$, of the normalized sum $\frac{1}{n}\sum_{i=1}^n
A_{p;F}(X_i)$ has been established by Avram and Taqqu \cite
{AvramTaqqu1987} for $1\le
p<1/(2\beta-1)$ (for the meaning of $\beta$ see Theorem \ref{nclt}
below). So, whenever the process ${\cal E}_{n,p;F}(\cdot)$ can be shown
to converge in probability to zero under the same re-scaling, we obtain
that the limit in distribution of a re-scaled version of the process
${\cal E}_{n,p-1;F}(\cdot)$ is given by
$(-1)^{p} F^{(p)}(\cdot) Z_{p,\beta}$. This idea is basically due to
Dehling and Taqqu \cite{DehlingTaqqu1989} who considered the Gaussian
case and the
uniform sup-metric $d_\infty$. For the linear process and the uniform
sup-metric $d_\infty$ this approach was used by Ho and Hsing \cite
{HoHsing1996} and
Wu \cite{Wu2003} for arbitrary $p\ge1$, and by
Giraitis and Surgailis \cite{GiraitisSurgailis1999} for $p=1$. Wu \cite
{Wu2003} also considered
bounds for the second moment of weighted sup-norms of the leading term
of ${\cal E}_{n,p-1;F}(\cdot)$. For the linear process and the weighted
sup-metric $d_{\phi_\lambda}$ the approach of Dehling and Taqqu \cite
{DehlingTaqqu1989}
was applied to the case $p=1$ by Beutner, Wu and Z{\"a}hle \cite
{BeutnerWuZaehle}. In the
following theorem, we generalize the latter to arbitrary $p\ge1$.

\begin{theorem}\label{nclt}
Let $p\in\N$, $\lambda\ge0$, and assume that:
\begin{itemize}[(a)]
\item[(a)] $a_s=s^{-\beta} \ell(s)$, $s\in\N$, where $\beta\in
(\frac{1}{2},1)$ and $\ell$ is slowly varying at infinity.
\item[(b)] $\ex[|\varepsilon_1|^{(4+2\lambda)\vee(2p)}]<\infty$.
\item[(c)] The df $G$ of $\varepsilon_1$ is $p+1$ times
differentiable and $\sum_{j=1}^{p+1}\int_\R|G^{(j)}(x)|^2\phi
_{2\lambda}(x) \mrmd x<\infty$.
%
\item[(d)] $p(2\beta-1)<1$.
\end{itemize}
Then
\[
\bigl\{n^{p(\beta-1/2)}\ell(n)^{-p}\bigr\} {\cal E}_{n,p-1;F}(
\cdot) \stackrel{{\mathsf d}} {\longrightarrow} (-1)^{p} F^{(p)}(
\cdot) Z_{p,\beta} \qquad\mbox{$\bigl($in $(\D_{\phi_\lambda},{\cal
D}_{\phi
_\lambda},d_{\phi_\lambda})\bigr)$},
\]
where
\[
Z_{p,\beta}:= c_{p,\beta}\int_{-\infty<u_1<\cdots<u_p<1} \Biggl\{
\int_0^1\eins_{(u_p,1)}(v)\prod
_{j=1}^p(v-u_j)^{-\beta}\mrmd v
\Biggr\} W(\dd u_1)\cdots W(\dd u_p)
\]
with $W$ a white noise measure (i.e., an additive Gaussian random set
function satisfying $\ex[W(B)]=0$ and $\ex[W(B)\cap W(B')]=|B\cap
B'|$ for all $B,B'\in{\cal B}(\R)$) and
\[
c_{p,\beta}:= \biggl(\frac{p! (1-p(\beta-{1/2}))
(1-p(2\beta-1))}{\int_0^\infty(x+x^2)^{-\beta}\mrmd x} \biggr)^{1/2}.
\]
\end{theorem}

\begin{remark}\label{VergleichDehlingTaqqu}
(i) The infinite moving average representation of an ARFIMA$(p,d,q)$
process with fractional difference parameter $d\in(0,1/2)$ satisfies
assumption (a) with $\beta=1-d$; see, for instance, Hosking \cite{Hosking1981},
Section 3.

\mbox{}\hphantom{i}(ii) Here we have chosen to define the stochastic process (\ref
{EqVarepsprocess}) in terms of $\hat{F}_n$, $F$ and the Appell
polynomials of $F$, because this allows us later to define the
statistic (\ref{eqcorlongmemory}) in terms of the df of the
observables. However, we conjecture that assumption (b) can be relaxed
to $\ex[|\varepsilon_1|^{(2+2\lambda)}]<\infty$ by replacing as in
Wu \cite{Wu2003,Wu2006}
the Appell\vspace*{1pt} polynomials $A_{j;F}(X_i)$ by the expressions
$A_{j;F}^{=(1,\ldots,1)}(X_i)$ (to be introduced at the beginning of the
proof of Theorem \ref{nclt}) and by the method of proof as in
Beutner, Wu and Z{\"a}hle \cite{BeutnerWuZaehle} where this was done
for $p=1$.


(iii) Condition (c) implies in particular that the df $F$ of $X_1$
is $p$ times differentiable with $F^{(p)}\in\D_{\phi_\lambda}$;
cf.
inequality (30) in Wu \cite{Wu2003} with $n=\infty$, $\kappa=1$ and
$\gamma=2\lambda$. Further, assumption (c) can be relaxed in that it
suffices to require that there is some $m\in\N$ such that the df
$G_m$ of $\sum_{s=0}^{m-1}a_s\varepsilon_{m-s}$ is $p+1$ times
differentiable and satisfies $\sum_{j=1}^{p+1}\int_\R|
G^{(j)}_m(x)|^2\phi_{2\lambda}(x) \mrmd x<\infty$. The proof still
works in this setting; see also Wu \cite{Wu2003}.
\end{remark}


\begin{pf*}{Proof of Theorem \ref{nclt}}
It was shown by Avram and Taqqu \cite{AvramTaqqu1987}, Theorem 1, that
the $p$th
Appell polynomial of $F$ evaluated at $X_i$ has the representation
\begin{eqnarray*}
A_{p;F}(X_i) & = & \sum_{\ell=1}^p
\sum_{q(\ell)\in\Pi_{\ell,p}}\frac
{p!}{q_1!\cdots q_{\ell}!} \sum
_{m(\ell)\in\Lambda_{q(\ell
)}}\prod_{k=1}^{\ell}
a_{m_k}^{q_k} A_{q_k;G}(\varepsilon_{i-m_k})
\\
& = & \sum_{m(p)\in\Lambda_{q(p)}} p! \prod
_{k=1}^p a_{m_k}A_{q_k;G}(
\varepsilon_{i-m_k})
\\
& &{} +\sum_{\ell=1}^{p-1} \sum
_{q(\ell)\in\Pi_{\ell,p}}\frac
{p!}{q_1!\cdots q_{\ell}!} \sum_{m(\ell)\in\Lambda_{q(\ell)}}
\prod_{k=1}^{\ell} a_{m_k}^{q_k}
A_{q_k;G}(\varepsilon_{i-m_k})
\\
&=:& A_{p;F}^{=(1,\ldots,1)}(X_i) + A_{p;F}^{\not=(1,\ldots,1)}(X_i),
\end{eqnarray*}
where $A_{q_k;G}$ denotes the $q_k$th Appell polynomial of the df
$G$ of $\varepsilon_1$, and, for every $\ell\in\{1,\ldots,p\}$,
$\Pi_{\ell,p}$ is the set of all
$q(\ell)=(q_1,\ldots,q_{\ell})\in\N^{\ell}$ satisfying
$q_1+\cdots+q_{\ell}=p$ and $1\le q_1\le\cdots\le q_{\ell}$.
Moreover, for a given $q(\ell)=(q_1,\ldots,q_{\ell})$ we denote by
$\Lambda_{q(\ell)}$ the set of all
$m(\ell)=(m_1,\ldots,m_{\ell})\in\N^{\ell}_0$ such that
$m_i\not=m_j$ for $i\not=j$ and, in addition, if $q_i=q_{i+1}$, then
$m_i<m_{i+1}$. So, introducing a telescoping sum, we obtain from
(\ref{generalizedempproc})
\begin{eqnarray*}
{\cal E}_{n,p-1;F}(\cdot) & = & \Biggl\{{\cal E}_{n,p;F}(\cdot)+
\sum_{j=1}^p (-1)^{j}
F^{(j)}(\cdot) \Biggl(\frac{1}{n}\sum
_{i=1}^n A_{j;F}^{\not
=(1,\ldots,1)}(X_i)
\Biggr) \Biggr\}
\\
& &{} - \sum_{j=1}^p (-1)^{j}
F^{(j)}(\cdot) \Biggl(\frac{1}{n}\sum
_{i=1}^n A_{j;F}^{\not=(1,\ldots,1)}(X_i)
\Biggr)
\\
& &{} + (-1)^{p} F^{(p)}(\cdot) \Biggl(\frac{1}{n}\sum
_{i=1}^n A_{p;F}(X_i)
\Biggr)
\\
& = & \Biggl\{\hat F_n(\cdot)-F(\cdot)- \sum
_{j=1}^p (-1)^{j} F^{(j)}(
\cdot) \Biggl(\frac{1}{n}\sum_{i=1}^n
A_{j;F}^{=(1,\ldots,1)}(X_i) \Biggr) \Biggr\}
\\
& &{} - \sum_{j=1}^{p-1} (-1)^{j}
F^{(j)}(\cdot) \Biggl(\frac
{1}{n}\sum
_{i=1}^n A_{j;F}^{\not=(1,\ldots,1)}(X_i)
\Biggr)
\\
& &{} + (-1)^{p} F^{(p)}(\cdot) \Biggl(\frac{1}{n}\sum
_{i=1}^n A_{p;F}^{=(1,\ldots,1)}(X_i)
\Biggr)
\\
& =: & S_{n,p}^0(\cdot) + T_{n,p}(\cdot) +
U_{n,p}(\cdot).
\end{eqnarray*}
Under assumptions (a), $\ex[|\varepsilon_0|^{2p}]<\infty$, and (d),
it follows from Step 3 in the proof of Theorem~2 of Avram and Taqqu
\cite{AvramTaqqu1987} that the expression
\[
\bigl\{n^{p(\beta-1/2)}\ell(n)^{-p}\bigr\} \frac{1}{n}\sum
_{i=1}^n A_{j;F}^{\not=(1,\ldots,1)}(X_i)
= \frac{1}{n^{1-p(\beta
-1/2)}\ell(n)^{p}}\sum_{i=1}^n
A_{j;F}^{\not=(1,\ldots,1)}(X_i)
\]
converges in probability to zero for every $j=1,\ldots, p$. So, in
view of $F^{(j)}\in\D_{\phi_\lambda}$, $j=1,\ldots,p$, we obtain
\[
\bigl\{n^{p(\beta-1/2)}\ell(n)^{-p}\bigr\} \mrmd _{\phi_\lambda}
\bigl(T_{n,p}(\cdot),0\bigr) \stackrel{\mathsf p} {\longrightarrow} 0.
\]
Avram and Taqqu \cite{AvramTaqqu1987}, Theorem 2, also showed that
under the same
assumptions the expression
\[
\bigl\{n^{p(\beta-1/2)}\ell(n)^{-p}\bigr\}\frac{1}{n}\sum
_{i=1}^n A_{p;F}^{=(1,\ldots,1)}(X_i)
= \frac{1}{n^{1-p(\beta-1/2)}\ell
(n)^{p}}\sum_{i=1}^n
A_{p;F}^{=(1,\ldots,1)}(X_i)
\]
converges in distribution to $Z_{p,\beta}$; for the shape of the
normalizing constant $c_{p,\beta}$ see Ho and Hsing \cite{HoHsing1996},
Lemma 6.1.
So, in view of $F^{(p)}\in\D_{\phi_\lambda}$, the process
$U_{n,p}(\cdot)$ converges in distribution to $(-1)^{p}
F^{(p)}(\cdot) Z_{p,\beta}$ w.r.t. $d_{\phi_\lambda}$.
In the remainder of the proof, we will show that
%
\begin{equation}
\label{nclt-proof-1} \bigl\{n^{p(\beta-1/2)}\ell(n)^{-p}\bigr\}
d _{\phi_\lambda}\bigl(S_{n,p}^0(\cdot),0\bigr) \stackrel{
\mathsf p} {\longrightarrow} 0
\end{equation}
so that assertion of Theorem \ref{nclt} will follow from Slutzky's lemma.

For (\ref{nclt-proof-1}) to be true, it suffices to show that
$d_{\phi_\lambda}(\frac{S_{n,p}(\cdot)}{\sigma_{n,p}},0)$
converges in
probability to zero, where $S_{n,p}:=n S_{n,p}^0$ and
$\sigma_{n,p}:=n^{1-p(\beta-1/2)}\ell(n)^{p}$. Under the assumptions
(a), $\ex[|\varepsilon_0|^{(4+2\lambda)}]<\infty$, (c) and (d), we have
from Theorem 2 and Lemma 5 of Wu \cite{Wu2003} that
%
\begin{equation}
\label{wulemma} \ex\bigl[d_{\phi_{2\lambda}}\bigl(S_{n,p}(
\cdot)^2,0\bigr) \bigr] = \RMO \bigl(n(\log n)^2+
\Xi_{n,p} \bigr)
\end{equation}
with $\Xi_{n,p}=\RMO(n^{2-(p+1)(2\beta-1)}\ell(n)^{2(p+1)})$
(notice that there is a typo in Lemma 5 of Wu \cite{Wu2003} where it must
be $p(2\beta-1)<1$ instead of $(p+1)(2\beta-1)<1$). From (\ref
{wulemma}), we obtain by the Markov inequality for some constant $C>0$ and
every $\varepsilon>0$
\begin{eqnarray*}
\pr\bigl[\sigma_{n,p}^{-1} d _{\phi_\lambda}
\bigl(S_{n,p}(\cdot),0\bigr)>\varepsilon\bigr] & = & \pr\bigl[
\sigma_{n,p}^{-2} d _{\phi_{2\lambda}}\bigl(S_{n,p}(
\cdot)^2,0\bigr)>\varepsilon^2 \bigr]
\\
& \le& \frac{1}{\varepsilon^2} \frac{\ex[d_{\phi_{2\lambda
}}(S_{n,p}(\cdot)^2,0) ]}{\sigma_{n,p}^2}
\\
& \le& C \varepsilon^{-2} \frac{n(\log n)^2+n^{2-(p+1)(2\beta
-1)}\ell(n)^{2(p+1)}}{n^{2-p(2\beta-1)}\ell(n)^{2p}}
\\
& \le& C \varepsilon^{-2} \biggl(\frac{(\log n)^2}{n^{1-p(2\beta
-1)}}+\frac{1}{n^{2\beta-1}}
\biggr)\ell(n)^2.
\end{eqnarray*}
Due to assumption (d), the latter bound converges to zero as $n\to
\infty$. That is, $d_{\phi_\lambda}(\frac{S_{n,p}(\cdot)}{\sigma
_{n,p}},0)$ indeed converges in probability to zero.
\end{pf*}

Combining Theorems \ref{theoremcontmappingforu-stat} and
\ref{nclt}, one can in principle easily derive the asymptotic
distribution of non-degenerate and degenerate $V$-statistics based on
linear long-memory sequences. For non-degenerate $V$-statistics (as,
e.g., Gini's mean difference from Example \ref{examplegini})
one can apply part (i) of Theorem \ref{theoremcontmappingforu-stat};
see also Hsing \cite{Hsing2000} who uses a different approach
for Gini's mean difference. For degenerate $V$-statistics (as, e.g., the
Cram\'{e}r--von Mises statistic from Example~\ref{gof}), one can apply
part (ii) of Theorem \ref{theoremcontmappingforu-stat}. However,
in the long-memory case the situation is often more complex because
several $V$-statistics based on long-memory sequences systematically
degenerate asymptotically. For instance, the (sample) variance with
corresponding kernel $g(x_1,x_2)=\frac{1}{2}(x_1-x_2)^2$ is
typically non-degenerate w.r.t. $(g,F)$ (cf. Example \ref
{examplevariance}), but in this case the integral on the right-hand
side in
(\ref{theoremcontmappingforu-stat-eq1}) with
$B^\circ(\cdot)=(-1)F^{(1)}(\cdot)Z_{1,\beta}$ equals
\[
-\sum_{i=1}^2\int
B^\circ(x-) \mrmd g_{i,F}(x) = Z_{1,\beta}\sum
_{i=1}^2\int F^{(1)}(x-) \bigl(x-
\ex[X_1]\bigr) \mrmd x = 0.
\]
Indeed, from Example \ref{examplevariance} we know that in this
case $\dd g_{1,F}(x)=\dd g_{2,F}(x)=(x-\ex[X_1]) \mrmd x$ holds, and hence
$\int
F^{(1)}(x-)(x-\ex[X_1])\mrmd x=\int F^{(1)}(x)(x-\ex[X_1])\mrmd x=0$. 
That is, in the long-memory case the sample variance regarded as a
$V$-statistic is asymptotically degenerate w.r.t.
$(g,F,(n^{(\beta-1/2)}\ell(n)^{-1})_n)$ in the sense of Section
\ref{secnotionofdegeneracyandnon-dgeneracy}, and so an application of
part (i) of Theorem~\ref{theoremcontmappingforu-stat} yields little.
Moreover, part (ii) of Theorem \ref {theoremcontmappingforu-stat} is
useful neither in this case, because part (ii) of Theorem
\ref{theoremcontmappingforu-stat} is based on the fact that the linear
part in the representation (\ref {EqBZdecomposition}) vanishes.
However, this is not the case here, since the sample variance is not
(finite sample) degenerate w.r.t. $(g,F)$. This is in accordance with
the remarkable observation of Dehling and Taqqu \cite{DehlingTaqqu1991}
that in the long-memory case both terms of the von Mises (resp.,
Hoeffding) decomposition of the sample variance contribute to the
asymptotic distribution. Dehling and Taqqu \cite{DehlingTaqqu1991}
considered the sample variance based on Gaussian long-memory sequences.
From the following Corollary~\ref{corollarytothmnclt-1}, we cannot only
derive the analogue for linear long-memory sequences (see Example \ref
{examplevariance-longmem}), but can also derive the asymptotic
distribution of more general asymptotically degenerate $U$- and
$V$-statistics based on linear long-memory sequences (see, e.g.,
Examples \ref{examplesquaredabsoluetvalue-longmem}, \ref
{exampleartificial-longmem} and \ref{exampletfs-longmem}). We note that
recently L{\'e}vy-Leduc \textit{et al.}  \cite{Levy-LeducAoS} also
derived the asymptotic distribution of some asymptotically degenerate
$U$-statistics (with bounded kernels) based on Gaussian long-memory
sequences using different techniques. For further applications of their
results, see L{\'e}vy-Leduc \textit{et al.}
\cite{Levy-LeducStatistics}.

\begin{corollary}\label{corollarytothmnclt-1}
Let $F$ be a df on the real line, and $g\dvtx \R^2\to\R$ be some
measurable function. Assume that the representation (\ref
{EqBZdecomposition}) with $F_n:=\hat F_n$ holds for $F$ and $g$, and that
%
\begin{equation}
\label{corollarytothmnclt-1-eq1} \sum_{i=1}^2
\int\phi_{-\lambda}(x) |\dd g_{i,F}|(x)<\infty
\quad\mbox{and}\quad \iint
\phi_{-\lambda}(x_1) \phi_{-\lambda
}(x_2)
|\dd g|(x_1,x_2)<\infty
\end{equation}
holds for some $\lambda\ge0$. Let $p,q,r\in\N$, set
%
\begin{eqnarray}
\label{eqcorlongmemory}
&&\hspace*{-22pt} {\cal V}_{n,g;p,q,r}(\hat{F}_n)
\nonumber\\
&&\hspace*{-22pt}\quad:= V_g(\hat F_n)-V_g(F) 
+ \sum
_{\ell=1}^2\sum_{j=1}^{p-1}(-1)^{j}
\Biggl(\frac{1}{n}\sum_{i=1}^nA_{j;F}(X_i)
\Biggr)\int F^{(j)}(x-) \mrmd g_{\ell,F}(x)
\nonumber
\\
&&\hspace*{-22pt}\qquad{} -\sum_{j=1}^{q-1}(-1)^{j}
\Biggl(\frac{1}{n}\sum_{i=1}^nA_{j;F}(X_i)
\Biggr) \iint F^{(j)}(x_1-) \bigl(\hat F_n(x_2-)-F(x_2-)
\bigr) \mrmd g(x_1,x_2)
\nonumber\\[-8pt]\\[-8pt]
&&\hspace*{-22pt}\qquad{} -\sum_{k=1}^{r-1}(-1)^{k}
\Biggl(\frac{1}{n}\sum_{i=1}^nA_{k;F}(X_i)
\Biggr) \iint\bigl(\hat F_n(x_1-)-F(x_1-)
\bigr) F^{(k)}(x_2-) \mrmd g(x_1,x_2)
\nonumber
\\
&&\hspace*{-22pt}\qquad{} +\sum_{j=1}^{q-1}\sum
_{k=1}^{r-1}(-1)^{j+k} \Biggl(
\frac{1}{n}\sum_{i=1}^nA_{j;F}(X_i)
\Biggr) \Biggl(\frac{1}{n}\sum_{i=1}^nA_{k;F}(X_i)
\Biggr)\nonumber\\
&&\hspace*{-22pt}\qquad\hspace*{39pt}{}\times\iint F^{(j)}(x_1-) F^{(k)}(x_2-)
\mrmd g(x_1,x_2)
\nonumber
\end{eqnarray}
(with the convention $\sum_{j=1}^0(\cdots):=0$), and assume that all
integrals on the right-hand side in (\ref{eqcorlongmemory}) are
well defined (which is in particular the case if $F$ is $\max\{p,q,r\}
$ times differentiable with $F^{(k)}\in\D_{\phi_\lambda}$ for all
$k=0,\ldots,\max\{p,q,r\}$).

\begin{longlist}
\item
Assume $q+r>p$ and that the assumptions \textup{(a)--(c)} of Theorem \ref
{nclt} with $p$ replaced by $\max\{p,q,r\}<\infty$ hold for the same
$\lambda$. Then, if in addition $s(2\beta-1)<1$ holds for each $s\in
\{p,q,r\}$, we have with $Z_{p,\beta}$ defined as in Theorem \ref{nclt}
\[
\bigl\{n^{p(\beta-1/2)}\ell(n)^{-p}\bigr\} {\cal V}_{n,g;p,q,r}(
\hat{F}_n)  \stackrel{\mathsf d} {\longrightarrow}  (-1)^{p}
Z_{p,\beta} \sum_{\ell=1}^2\int
F^{(p)}(x-) \mrmd g_{\ell,F}(x). 
\]

\item Assume $q+r=p$ and that the assumptions \textup{(a)--(d)} of Theorem \ref
{nclt} hold for the same $\lambda$ and $p$. Then we have with
$Z_{s,\beta}$ defined as in Theorem \ref{nclt} for $s\in\{p,q,r\}$
%
\begin{eqnarray}
\label{eqcortoncltequalp,q,r}\label{corollarytothmnclt-1-eq3}
&&
\bigl\{n^{p(\beta-1/2)}\ell(n)^{-p}\bigr
\}{\cal V}_{n,g;p,q,r}(\hat{F}_n) \nonumber\\
&&\quad \stackrel{\mathsf d} {
\longrightarrow}  (-1)^{p} Z_{p,\beta} \sum
_{\ell=1}^2\int F^{(p)}(x-)
\mrmd g_{\ell,F}(x)
\\
&&\hphantom{\stackrel{\mathsf d} {
\longrightarrow}}\quad{}+ (-1)^{p} Z_{q,\beta}Z_{r,\beta}\iint
F^{(q)}(x_1-)F^{(r)}(x_2-)
\mrmd g(x_1,x_2).
\nonumber
\end{eqnarray}
\end{longlist}
\end{corollary}

Recall that assumption (c) of Theorem \ref{nclt} implies that $F$ is
$\max\{p,q,r\}$ times differentiable and that all derivatives up to
the $\max\{p,q,r\}$th derivative lie in $\D_{\phi_\lambda}$.

\begin{remark}
The random variables $Z_{p,\beta}$, $Z_{q,\beta}$ and $Z_{r,\beta}$ in
part (ii) of Corollary \ref{corollarytothmnclt-1} are dependent.
The specification of their joint distribution seems to be an open
problem. Only in the Gaussian case the joint cumulants of
$Z_{1,\beta}^2$ and $Z_{2,\beta}$ are known from the supplementary
material to L{\'e}vy-Leduc \textit{et al.}  \cite{Levy-LeducAoS}. Notice that it
is even hard to
specify the (Rosenblatt) distribution of~$Z_{2,\beta}$; for details see
Veillette and Taqqu \cite{TaqquVeillette2011}.
\end{remark}

\begin{pf*}{Proof of Corollary \ref{corollarytothmnclt-1}}
Using the representation (\ref{EqBZdecomposition}) of $V_g(\hat
F_n)-V_g(F)$, we obtain that ${\cal V}_{g,n;p,q,r}(F)$ equals
\begin{eqnarray*}
&& -\sum_{i=1}^2\int\bigl[\hat
F_n(x-)-F(x-)\bigr] \mrmd g_{i,F}(x) +\iint(\hat
F_n-F) (x_1-) (F_n-F) (x_2-)
\mrmd g(x_1,x_2)
\\
&&\qquad\hspace*{0pt}{} +\sum_{\ell=1}^2\sum
_{j=1}^{p-1}(-1)^{j} \Biggl(
\frac{1}{n}\sum_{i=1}^nA_{j;F}(X_i)
\Biggr)\int F^{(j)}(x-) \mrmd g_{\ell,F}(x)
\\
&&\qquad\hspace*{0pt}{} -\sum_{j=1}^{q-1}(-1)^{j}
\Biggl(\frac{1}{n}\sum_{i=1}^nA_{j;F}(X_i)
\Biggr)\iint F^{(j)}(x_1-) \bigl(\hat F_n(x_2-)-F(x_2-)
\bigr) \mrmd g(x_1,x_2)
\\
&&\qquad\hspace*{0pt}{} -\sum_{k=1}^{r-1}(-1)^{k}
\Biggl(\frac{1}{n}\sum_{i=1}^nA_{k;F}(X_i)
\Biggr)\iint\bigl(\hat F_n(x_1-)-F(x_1-)
\bigr) F^{(k)}(x_2-) \mrmd g(x_1,x_2)
\\
&&\qquad\hspace*{0pt}{} +\sum_{j=1}^{q-1}\sum
_{k=1}^{r-1}(-1)^{j+k} \Biggl(
\frac{1}{n}\sum_{i=1}^nA_{j;F}(X_i)
\Biggr) \Biggl(\frac{1}{n}\sum_{i=1}^nA_{k;F}(X_i)
\Biggr)\\
&&\qquad\hspace*{39pt}{}\times\iint F^{(j)}(x_1-) F^{(k)}(x_2-)
\mrmd g(x_1,x_2)
\\
&&\quad\hspace*{0pt}{} = -\sum_{\ell=1}^2\int{\cal
E}_{n,p-1;F}(x-) \mrmd g_{\ell,F}(x) + \iint{\cal E}_{n,q-1;F}(x_1-)
{\cal E}_{n,r-1;F}(x_2-) \mrmd g(x_1,x_2).
\end{eqnarray*}
Moreover, by Theorem \ref{nclt},
\[
\bigl\{n^{s(\beta-1/2)}\ell(n)^{-s}\bigr\}{\cal E}_{n,s-1;F}(
\cdot) \stackrel{\mathsf d} {\longrightarrow} (-1)^{s}c_{s,\beta}Z_{s,\beta
}F^{(s)}(
\cdot) \qquad\mbox{$\bigl($in $(\D_{\phi_\lambda},{\cal D}_{\phi
_\lambda},d_{\phi_\lambda})
\bigr)$}
\]
for $s=p,q,r$. Therefore, assertion (i) follows from the Continuous
Mapping theorem and (\ref{corollarytothmnclt-1-eq1}) as well
as Slutzky's lemma and the assumption $q+r>p$. Moreover, assertion (ii)
follows from the Continuous Mapping theorem, (\ref
{corollarytothmnclt-1-eq1}) and the assumption $p=q+r$.
\end{pf*}

It is worth pointing out that, as mentioned at the beginning of this
section, ${\cal V}_{n,g;p,q,r}(\hat{F}_n)$ is obtained by using the
representation (\ref{EqBZdecomposition}) and an ``expansion'' of
$\hat{F}_n-F$ in the sense of (\ref{EqVarepsprocess}). Obviously, with
increasing $p$, $q$ or $r$, the expression ${\cal V}_{n,g;p,q,r}(F)$
defined in (\ref{eqcorlongmemory}) is getting more and more involved.
So, for statistical applications one should choose $p$, $q$ and $r$ as
small as possible. On the other hand, $p$ has to be chosen so large so
that the limit in distribution of $\{
n^{p(\beta-1/2)}\ell(n)^{-p}\}{\cal V}_{n,g;p,q,r}(F)$ does not vanish.
That is, an application of Corollary \ref{corollarytothmnclt-1}
requires a trade-off between the simplicity of the statistic ${\cal
V}_{n,g;p,q,r}(F)$ and the benefit of the asymptotic distribution. A
particularly favorable situation is the one where some (or preferably
all) terms on the right-hand side of (\ref{eqcorlongmemory}), which
are different from $V_g(\hat F_n)-V_g(F)$, vanish. This is the case if
the respective integrals $\int F^{(j)}(x-) \mrmd g_{\ell,F}(x)$ etc. vanish,
for instance, in the case of the sample variance and in the case of the
test for symmetry; cf. Examples \ref{examplevariance-longmem} and
\ref{exampletfs-longmem}. In other situations, the statistic ${\cal
V}_{n,g;p,q,r}(F)$ might be more
complicated than $V_g(\hat F_n)-V_g(F)$. 
Yet, it seems to be among the best achievable results. Finally, notice
that in cases where part (i) or part (ii) of Theorem \ref
{theoremcontmappingforu-stat} already yields a non-trivial asymptotic
distribution, the result can also be derived by Corollary \ref
{corollarytothmnclt-1}. This is exemplified in the next remark.

\begin{remark}
(i) For Gini's mean difference take $p=q=r=1$. Then the result obtained
from part (i) of Corollary \ref{corollarytothmnclt-1} coincides
with the result we get from part (i) of Theorem~\ref
{theoremcontmappingforu-stat}.

(ii) For the weighted Cram\'{e}r--von Mises statistic take $p=2$ and
$q=r=1$. Then, under the hypothesis that $F=F_0$, we have from Corollary
\ref{corollarytothmnclt-1}(ii) that the asymptotic distribution
equals $Z_{1,\beta}Z_{1,\beta}\iint F^{(1)}(x_1-)F^{(1)}(x_2-)
\mrmd g(x_1,x_2)= (Z_{1,\beta})^2\int w(x)(F^{(1)}(x-))^2 \mrmd F(x)$,
where we used (\ref{tfs-eq-2}). That is, in accordance with Example
\ref{exampleGiniandCramer-vonMisesdegenerate}.
\end{remark}

Gini's mean difference discussed in part (i) of the preceding remark
is an example for an asymptotically non-degenerate $U$- or
$V$-statistic. The weighted Cram\'{e}r--von Mises statistic discussed in
part (ii) of the preceding remark is an example for an
asymptotically degenerate $U$- or $V$-statistic of type (1.a) in the
sense of Section \ref{secnotionofdegeneracyandnon-dgeneracy}.
The following two Examples \ref{examplesquaredabsoluetvalue-longmem}
and \ref{examplevariance-longmem} provide
some asymptotically degenerate $U$- or $V$-statistics of type (1.b) and
type (1.c), respectively. Examples \ref{exampleartificial-longmem}
and \ref{exampletfs-longmem}
below will provide some asymptotically degenerate $U$-
or $V$-statistics of type 2 in the sense of Section \ref
{secnotionofdegeneracyandnon-dgeneracy}.

\begin{example}[(Squared absolute mean of a symmetric
distribution)]\label{examplesquaredabsoluetvalue-longmem}
The kernel $g(x_1,x_2)=x_1\cdot x_2$ for estimating the
squared mean has been investigated repeatedly in the literature. Let
us consider here the related kernel $g(x_1,x_2)=|x_1|\cdot|x_2|$ for
estimating the squared \textit{absolute} mean of a distribution $F$
having a finite first moment. In this case, we obtain
$g_{i,F}(x)=\ex[|X_1|]\cdot|x|$ and hence $\dd g_{i,F}(x) =
\ex[|X_1|](-\eins_{\{x<0\}}+\eins_{\{x\ge0\}}) \mrmd x$
for $i=1,2$. Moreover, we have $\dd g(x_1,x_2)=(
\eins_{\{x_1\geq0, x_2\geq0\}}-\eins_{\{x_1<0,x_2 \geq
0\}}-\eins_{\{x_1\geq0, x_2<0\}}+\eins_{\{x_1<0, x_2<0\}})
\mrmd x_1\mrmd x_2$. It is then
easily checked that the conditions of Lemmas \ref
{Lemmaintegrationbypartsonedimension} and \ref
{Lemmaintegrationbypartstwodimensions} are fulfilled for any weight
function $\phi$ with $\int
1 \slash\phi(x) \mrmd x < \infty$. Now, let us in addition assume that
$F^{(1)}$ is symmetric about $0$. Then, on one hand,
Theorem \ref{theoremcontmappingforu-stat}(i) with
$B^\circ(\cdot)=(-1)F^{(1)}(\cdot)Z_{1,\beta}$ yields that $\{
n^{\beta-1/2}\ell(n)^{-1}\}(V_g(\hat
F_n)-V_g(F))$ converges in distribution to
$-\sum_{i=1}^2\int B^\circ(x-) \mrmd g_{i,F}(x)=2Z_{1,\beta}\ex
[|X_1|](-\int_{-\infty}^0 F^{(1)}(x-) \mrmd x +\int_0^{\infty}
F^{(1)}(x-) \mrmd x=0$.
On the other hand, part (ii) of Theorem \ref
{theoremcontmappingforu-stat} is helpful neither because it only yields
that $\{
n^{2(\beta-1/2)}\ell(n)^{-2}\}(V_g(\hat F_n)-V_g(F))$ converges in
distribution to
$\iint B^\circ(x_1-)B^\circ(x_2-) \mrmd g(x_1,x_2)=Z_{1,\beta}^2\iint
F^{(1)}(x_1)F^{(1)}(x_2) \mrmd g(x_1,x_2)=0$,
where for the latter ``$=$'' we used the symmetry of $F^{(1)}$.
However, if we take $p=2$ and $q=r=1$, we have ${\cal
V}_{n,g;2,1,1}(\hat{F}_n)=V_g(\hat F_n)-V_g(F)$ and obtain by
Corollary \ref{corollarytothmnclt-1}(ii)
\begin{eqnarray*}
&&\bigl\{n^{2\beta-1}\ell(n)^{-2}\bigr\} \bigl(V_g(\hat
F_n)-V_g(F) \bigr)
\\
&&\quad \stackrel{\mathsf d} {\longrightarrow}  Z_{2,\beta} \sum
_{\ell
=1}^2\int F^{(2)}(x-)
\mrmd g_{\ell,F}(x) + Z_{1,\beta}^2\iint F^{(1)}(x_1-)F^{(1)}(x_2-)
\mrmd g(x_1,x_2)
\\
&&\quad =  2 Z_{2,\beta} \biggl(-\int_{-\infty}^0
F^{(2)}(x) \mrmd x + \int_0^{\infty}F^{(2)}(x)
\mrmd x \biggr) = 4 Z_{2,\beta} \biggl(\int_0^{\infty}F^{(2)}(x)
\mrmd x \biggr),
\end{eqnarray*}
where for the latter ``$=$'' we used the antisymmetry of $F^{(2)}$
(i.e., $F^{(2)}(x)=-F^{(2)}(-x)$) which holds by the symmetry of
$F^{(1)}$. This shows that in the present
case $V_g(\hat F_n)$ is an asymptotically degenerate $V$-statistic
w.r.t. $(g,F,(n^{(\beta-1/2)}\ell(n)^{-1})))$ of type (1.b) in the
sense of Section \ref{secnotionofdegeneracyandnon-dgeneracy}.
\end{example}

\begin{example}[(Variance)]\label{examplevariance-longmem}
As discussed above, in our long-memory setting we can neither apply
part (i) nor part (ii) of Theorem \ref{theoremcontmappingforu-stat} to
derive a \textit{non-trivial} asymptotic distribution for the
sample variance; recall that the sample variance is a $V$-statistic
with corresponding kernel $g(x_1,x_2)=\frac{1}{2}(x_1-x_2)^2$.
However, part (ii) of Corollary \ref{corollarytothmnclt-1}
enables us to derive a non-trivial asymptotic distribution. From
Example \ref{examplevariance}, we know that
$\dd g_{1,F}(x)=\dd g_{2,F}(x)=(x-\ex[X_1]) \mrmd x$, which implies $\int
F^{(1)}(x-) \mrmd g_{\ell,F}(x)=0$. So we have ${\cal
V}_{n,g;2,1,1}(F)=V_g(\hat F_n)-V_g(F)$ and obtain by Corollary
\ref{corollarytothmnclt-1}(ii)
%
\begin{eqnarray}\label{eqexamplevarlongmemory}
&&
\bigl\{n^{2\beta-1}\ell(n)^{-2}\bigr\} \bigl(V_g(\hat
F_n)-V_g(F) \bigr)
\nonumber
\\
&&\quad \stackrel{\mathsf d} {\longrightarrow}  Z_{2,\beta} \sum
_{\ell
=1}^2\int F^{(2)}(x-)
\mrmd g_{\ell,F}(x) + Z_{1,\beta}^2\iint F^{(1)}(x_1-)F^{(1)}(x_2-)
\mrmd g(x_1,x_2)\qquad
\nonumber\\[-8pt]\\[-8pt]
&&\quad =  2 Z_{2,\beta}\int F^{(2)}(x-) \bigl(x-\ex[X_1]
\bigr) \mrmd x - \biggl(Z_{1,\beta}\int F^{(1)}(x) \mrmd x
\biggr)^2
\nonumber
\\
&&\quad = 2 Z_{2,\beta}\int F^{(2)}(x-) \bigl(x-\ex[X_1]
\bigr) \mrmd x - Z_{1,\beta
}^2,\nonumber
\end{eqnarray}
where for the first ``$=$'' we used the fact that $\dd g(x_1,x_2)$ is
the negative of the Lebesgue measure on~$\R^2$; cf. Example
\ref{examplevariance}. Notice that $\int
F^{(1)}(x-) \mrmd g_{\ell,F}(x)=0$ holds for every (sufficiently smooth)
df $F$, so that (\ref{eqexamplevarlongmemory}) is fully
satisfactory even in a non-parametric setting. Notice also that in the
present case $V_g(\hat F_n)$ is an asymptotically degenerate
$V$-statistic w.r.t. $(g,F,(n^{(\beta-1/2)}\ell(n)^{-1}))$ of type
(1.c) in the sense of Section \ref
{secnotionofdegeneracyandnon-dgeneracy}.
\end{example}

\begin{example}[(Scaling sequence $\bolds{(a_n^3)}$)]\label
{exampleartificial-longmem}
Let us consider the kernel $g(x_1,x_2)=x_1(|x_2|-1)$, and suppose
that $F^{(1)}$ is symmetric about zero and that
$m:=\ex[|X_1|]=1$. This setting is somewhat artificial but it leads
to quite an interesting limiting behavior of the corresponding $U$- or
$V$-statistic. We have $g_{1,F}(x_1)=x_1(m-1)=0$ and
$g_{2,F}(x_2)=\ex[X_1](|x_2|-1)=0$ due to the assumption $m=1$ and
the symmetry of $F^{(1)}$, respectively. That is,
$V_g(\hat F_n)$ is a degenerate $V$-statistic w.r.t. $(g,F)$, and
consequently an application of part (i) of Theorem \ref
{theoremcontmappingforu-stat} does not lead to a non-degenerate
limiting distribution.
Moreover, it is easily seen that $\dd g(x_1,x_2)=(\eins_{\{x_2 \geq0\}
}-\eins_{\{x_2<0\}}) \mrmd x_1\mrmd x_2$.
Therefore, part (ii) of Theorem \ref{theoremcontmappingforu-stat} does
not provide a tool to derive a non-degenerate limiting
distribution,
because under the imposed assumption the right-hand side in (\ref
{theoremcontmappingforu-stat-eq2}) with $B^\circ(\cdot
)=(-1)F^{(1)}(\cdot)Z_{1,\beta}$
vanishes.

On the other hand, part (ii) of Corollary \ref{corollarytothmnclt-1}
enables us to derive a \textit{non-trivial} asymptotic distribution. In
contrast to the sample variance in Example \ref
{examplevariance-longmem}, however, it does not make sense to work with
${\cal
V}_{n,g;2,1,1}(F)$, because in this case the limit in (\ref
{corollarytothmnclt-1-eq3}) vanishes. Indeed, the first summand of the
limit vanishes since $g_{1,F}\equiv g_{2,F}\equiv0$, and the second
summand of the limit vanishes since under the imposed assumptions the
integral $\iint F^{(1)}(x_1-)F^{(1)}(x_2-) \mrmd g(x_1,x_2)$ equals zero.
As a consequence we need to work with a $p$ larger than $2$. For
instance, for $(p,q,r)=(3,1,2)$ we obtain from Corollary \ref
{corollarytothmnclt-1}(ii) and $g_{1,F}\equiv g_{2,F}\equiv0$ that
\begin{eqnarray*}
&&\bigl\{n^{3(\beta-1/2)}\ell(n)^{-3}\bigr\}{\cal V}_{n,g;3,1,2}(F)
\\
&&\quad \stackrel{\mathsf d} {\longrightarrow}  - Z_{3,\beta} \sum
_{\ell
=1}^2\int F^{(3)}(x-)
\mrmd g_{\ell,F}(x) - Z_{1,\beta}Z_{2,\beta
}\iint
F^{(1)}(x_1-)F^{(2)}(x_2-)
\mrmd g(x_1,x_2)
\\
&&\quad = - Z_{1,\beta}Z_{2,\beta} \biggl(\int_0^{\infty}
\int F^{(1)}(x_1)F^{(2)}(x_2)
\mrmd x_1\mrmd x_2 -\int_{-\infty}^0
\int F^{(1)}(x_1)F^{(2)}(x_2)
\mrmd x_1\mrmd x_2 \biggr)
\\
&&\quad = -2 Z_{1,\beta}Z_{2,\beta} \int_0^\infty
F^{(2)}(x_2) \mrmd x_2,
\end{eqnarray*}
which is typically distinct from zero.
Notice that above we may replace ${\cal V}_{n,g;3,1,2}(F)$ by $V_g(\hat
F_n)-V_g(F)$ since
\begin{eqnarray*}
&&
{\cal V}_{n,g;3,1,2}(F)
\\
&&\quad =  V_g(\hat F_n)-V_g(F) + \sum
_{\ell=1}^2\sum_{j=1}^{2}(-1)^{j}
\Biggl(\frac{1}{n}\sum_{i=1}^nA_{j;F}(X_i)
\Biggr)\int F^{(j)}(x-) \mrmd g_{\ell,F}(x)
\\
&&\qquad{} + \Biggl(\frac{1}{n}\sum_{i=1}^nA_{1;F}(X_i)
\Biggr)\iint\bigl(\hat F_n(x_1-)-F(x_1-)
\bigr) F^{(1)}(x_2-) \mrmd g(x_1,x_2)
\\
&&\quad = V_g(\hat F_n)-V_g(F)+ \Biggl(
\frac{1}{n}\sum_{i=1}^nA_{1;F}(X_i)
\Biggr) \biggl(\int_0^{\infty}\int\bigl(\hat
F_n(x_1-)-F(x_1-)\bigr) F^{(1)}(x_2)
\mrmd x_1\mrmd x_2
\\
&&\qquad{} -\int_{-\infty}^0 \int\bigl(\hat
F_n(x_1-)-F(x_1-)\bigr) F^{(1)}(x_2)
\mrmd x_1\mrmd x_2 \biggr)
\\
&&\quad = V_g(\hat F_n)-V_g(F),
\end{eqnarray*}
where we used $g_{1,F}\equiv g_{2,F}\equiv0$, the continuity of
$F^{(1)}$, and the symmetry of $F^{(1)}$ about zero.
Thus, in the present case $V_g(\hat F_n)$ is an asymptotically degenerate
$V$-statistic w.r.t. $(g,F,(n^{(\beta-1/2)}\ell(n)^{-1}))$ of type 2
in the sense of Section \ref{secnotionofdegeneracyandnon-dgeneracy}.
\end{example}

The next example shows that it might even not be
sufficient to take the scaling sequence $(n^{3(\beta-1/2)})$ to obtain
a non-degenerate
limiting distribution.

\begin{example}[(Test for symmetry, scaling sequence
$\bolds{(a_n^4)}$)]\label{exampletfs-longmem}
Let us come back to the test statistic $T_n$ introduced in Example \ref
{tfs}, which is a $V$-statistic with kernel given by (\ref{tfsym-eq-1}).
We restrict to the null hypothesis that
the distribution is symmetric about zero. We have seen in Example
\ref{tfs} that in this case we obtain $g_{1,F}\equiv g_{2,F}\equiv
0$ and $\dd g(x_1,x_2)={\cal H}_{D}^1(d(x_1,x_2))-{\cal H}_{\widetilde
D}^1(d(x_1,x_2))$. That is, under the null hypothesis, $T_n$ can be
seen as a degenerate $V$-statistic. So, in principle, we could apply
Theorem \ref{theoremcontmappingforu-stat}(ii) to derive the
asymptotic distribution of $T_n=V_{g}(\hat F_n)$. However, the
integral on the right-hand side in (\ref
{theoremcontmappingforu-stat-eq2}) with $B^\circ(\cdot
)=(-1)F^{(1)}(\cdot)Z_{1,\beta}$
equals $\iint B^{\circ}(x_1)B^{\circ}(x_2) \mrmd g(x_1,x_2)=Z_{1,\beta
}^2(\int F^{(1)}(x)F^{(1)}(x) \mrmd x-\int F^{(1)}(x)F^{(1)}(-x) \mrmd x)=0$,
because $F^{(1)}$ is symmetric about zero. Now one might tend to apply
Corollary \ref{corollarytothmnclt-1} as in Example \ref
{examplevariance-longmem}, that is, with $(p,q,r)=(3,1,2)$, to obtain a
non-trivial limiting distribution. However, the integrals on the
right-hand side of (\ref{eqcortoncltequalp,q,r}) equal zero in
that case. Indeed, the first one equals zero, because $g_{1,F}\equiv
g_{2,F}\equiv
0$. The second one, which is given by $\int F^{(1)}(x)F^{(2)}(x)
\mrmd x-\int F^{(1)}(x)F^{(2)}(-x) \mrmd x$, equals zero, because of the
symmetry of $F^{(1)}$ and the antisymmetry of~$F^{(2)}$. However,
applying part (ii) of Corollary \ref{corollarytothmnclt-1} with
$(p,q,r)=(4,2,2)$ we obtain (using $g_{1,F}\equiv g_{2,F}\equiv
0$)
\begin{eqnarray*}
n^{4\beta-2}\ell(n)^{-4}{\cal V}_{n,g;4,2,2}(
\hat{F}_n) & \stackrel{\mathsf d} {\longrightarrow} &
Z_{2,\beta}^2 \biggl(\int F^{(2)}(x)F^{(2)}(x)
\mrmd x-\int F^{(2)}(x)F^{(2)}(-x) \mrmd x \biggr)
\\
& = & 4 Z_{2,\beta}^2 \int_0^\infty
\bigl(F^{(2)}(x) \bigr)^2 \mrmd x
\end{eqnarray*}
by the anti-symmetry of $F^{(2)}$. Using the symmetry of $F^{(1)}$ and
once again that $g_{1,F}\equiv g_{2,F}\equiv
0$, it can be easily checked that ${\cal V}_{n,g;4,2,2}(\hat
{F}_n)=V_g(\hat{F}_n)-V_g(F)$.
\end{example}


\section*{Acknowledgement}

The authors are very grateful to a referee for providing Example
\ref{examplereferee}, for his/her very careful reading, and
pointing out a missing word in Example \ref
{exampledegneratevariaceandsquared} that resulted in an incorrect statement.

\begin{supplement}
\stitle{Supplement to paper ``Continuous mapping approach to the
asymptotics of $U$- and $V$-statistics''}
\slink[doi]{10.3150/13-BEJ508SUPP} 
\sdatatype{.pdf}
\sfilename{BEJ508\_supp.pdf}
\sdescription{The supplement Beutner and Z{\"a}hle \cite{Beutnersuppl}
contains a discussion
of some extensions and limitations of the approach presented in this
paper.}
\end{supplement}

%

\printhistory

\end{document}